\documentclass[preprint,12pt]{elsarticle}

\usepackage[pdftex]{hyperref}
\biboptions{sort&compress,numbers,square}

\usepackage{amsfonts}
\usepackage{graphicx}
\usepackage{epsfig}
\usepackage{geometry}
\usepackage{tikz}
\usetikzlibrary{trees,positioning,shapes}
\usepackage{subcaption}
\usepackage{amsmath}
\usepackage{amsthm}
\usepackage{tabularx}
\usepackage{algorithmic}

\usepackage{xcolor, colortbl}

\usepackage{float}
\usepackage[hypcap]{caption}
\usepackage{subcaption}
\usepackage{algorithm2e}
\usepackage{array}
\usepackage{bm}

\usepackage{cancel}

\usepackage{siunitx}

\usepackage{url}

\usepackage{pdfpages}

\definecolor{cadmiumgreen}{rgb}{0.0, 0.42, 0.24}

\newcommand{\revision}[1]{\textcolor{black}{#1}}

\newcommand{\Par}[1]{\left({#1}\right)}

\newcommand{\Prod}[2]{\Par{#1, #2}}

\newcommand{\Dt}{\partial_t}

\newcommand{\D}[1]{\partial_{#1}}

\definecolor{medgreen}{rgb}{0,0.7,0}

\newcommand{\Sig}{\bm{\sigma}}

\DeclareMathOperator*{\DivX}{div_x}
\DeclareMathOperator*{\DivF}{div_{x,t}}

\newcommand{\BTau}{\boldsymbol{\tau}}
\newcommand{\SigU}{\boldsymbol{\underline{\sigma}}}
\newcommand{\BTauU}{\boldsymbol{\underline{\tau}}}
\newcommand{\Lop}{\mathcal{L}}

\usepackage{rotating}

\begin{document}

\begin{frontmatter}

\title{Solver algorithm for stabilized space-time formulation \\ of advection-dominated diffusion problem}

\author{Marcin \L{}o\'{s}$^1$, Paulina Sepulveda-Salas$^2$, Maciej Paszy\'nski$^1$}

\address{$^1$AGH University of Science and Technology, Krak\'ow, Poland \\
$^2$Pontificia Universidad Cat\'olica de Valparaiso, Chile}

\begin{abstract}
This article shows how to develop an efficient solver for a stabilized numerical space-time formulation of the advection-dominated diffusion transient equation. At the discrete space-time level, we approximate the solution by using higher-order continuous B-spline basis functions in its spatial and temporal dimensions. This problem is very difficult to solve numerically using the standard Galerkin finite element method due to artificial oscillations present when the advection term dominates the diffusion term. However, a first-order constraint least-square formulation allows us to obtain numerical solutions avoiding oscillations. The advantages of space-time formulations are the use of high-order methods and the feasibility of developing space-time mesh adaptive techniques on well-defined discrete problems.
We develop a solver for a least-square formulation to obtain a stabilized and symmetric problem on finite element meshes.
The computational cost of our solver is bounded by the cost of the inversion of the space-time mass and stiffness (with one value fixed at a point) matrices and the cost of the GMRES solver applied for the symmetric and positive definite problem.
We illustrate our findings on an advection-dominated diffusion space-time model problem and present two numerical examples: one with isogeometric analysis discretizations and the second one with an adaptive space-time finite element method.
\end{abstract}


\end{frontmatter}

\section{Introduction}
\label{sec:method}

The problem of developing an efficient solver for space-time formulations is important for many applications. The number of papers and citations in the field of space-time formulations is growing exponentially (Web of Science search ,,space-time formulation''), with some novel examples \cite{st1,st2,st3,st4,st5,st6,st7,st8,st9}. The direct solvers for space-time formulations are costly \cite{ST}. However, the use of space-time adaptive strategies can lower that cost \cite{st6}. In this paper, we study a space-time  formulation of an advection-dominated diffusion problem and derive a hybrid direct-iterative solver based on symmetric and positive definite matrices obtained from its discretization. 
The problem to be studied in this work is the following transient advection-dominated diffusion equation, defined on a domain $\Omega$  with boundary $\partial\Omega$ in a time interval $[0,T]$,
\begin{equation}
\begin{aligned}
  \Dt \phi -\varepsilon\Delta \phi +\beta \cdot \nabla \phi &= f(\mathbf{x},t), &(\mathbf{x},t)\in \Omega\times(0,T),  \\
  \phi(\mathbf{x}, 0) &= \phi_0(\mathbf{x} ) & \mathbf{x} \in \Omega,\\
  \phi(\mathbf{x}, t) &= 0 & (\mathbf{x}, t) \in \partial \Omega \times [0, T].
\end{aligned} \label{eq:problem1}
\end{equation}
This problem has important applications, from pollution simulations to modeling of flow and transport. Usually, it is solved numerically using the finite element method and time-stepping techniques. However, for small values of ~$\varepsilon / \left\|\beta\right\|$, which means the advection is much larger than diffusion, this problem becomes difficult to solve. Some finite element discretizations of this second-order problem do not work as desired since they encounter numerical instabilities resulting in unexpected oscillations and giving unphysical solutions.
Moreover, the traditional way of solving transient problems involves the introduction of a time-stepping scheme and the stabilization of solutions at every time step. In other words, the problem is translated into a sequence of problems representing the solutions at a given time interval: each problem in a time slab needs to be stabilized to solve the problem in the next time step. However, these formulations have certain limitations in terms of exploiting adaptive space-time strategies and parallelization of solvers due to their sequential nature.

There are several strategies to obtain stabilized methods: residual minimization \cite{c6}, streamline upwind Petrov-Galerkin \cite{c4}, and Discontinuous Petrov-Galerkin \cite{c5} methods. Some of them are developed for the stationary problem:
\begin{equation}
\begin{aligned}
  -\varepsilon\Delta \phi + \boldsymbol{\beta} \cdot \nabla \phi & = f(\mathbf{x}), & \mathbf{x}\in \Omega,  \\
  \phi(\mathbf{x}) &= 0 & \mathbf{x} \in \partial \Omega.
\end{aligned}
\end{equation}
A standard finite element method for this problem requires transforming this problem into a well-posed variational formulation, where we seek $u \in U\equiv V$ (where $U$ is an infinite dimensional space) such that
\begin{align*}
b(\phi,\psi)=l(\psi), \quad \forall \psi \in V, 
\end{align*}
where the bilinear form $b:U\times V\to \mathbb{R}$ and the linear form $l:V\to \mathbb{R}$ are such that the problem is well-posed. As an example, a finite element discretization can be obtained by defining the bilinear form $b$ and the linear form $f$ by 
\begin{align*}
b(\phi,\psi):= \int_{\Omega} (\boldsymbol{\beta} \cdot \nabla \phi\ \psi +\epsilon \nabla \phi \cdot \nabla \psi ) d{\bf x}, \quad\text{ and } \quad
l(\psi) := \int_{\Omega}
 f\psi d{\bf x}.
\end{align*}
The discrete form of this problem is obtained after selecting finite-dimensional spaces $U_h\subset U$ and $V_h\subset V$, and seeking for an approximated solution $u_h \in U_h \subset U$ such that
\begin{align*}
b(\phi_h,\psi_h)=l(\psi_h), \quad \forall \psi_h \in V_h \subset V.
\end{align*}
The study of the stability of a numerical method at the discrete level is also important to obtain stable numerical solutions for different material parameters.  The most important mathematical theorem allowing the study of the stability of the finite element method was independently proposed by prof. Ivo Babu\'ska, prof. Franco Brezzi, and prof. Olga \L{}ady\.zenskaya \cite{c1,c2,c3}.
They discovered (at about the same time independently of each other) equivalent conditions, now called inf-sup ("inf-sup condition"). The inf-sup condition is still used to this day by scientists who study the stability of the finite element method. This condition 
$\inf_{u\in  U} \sup_{v \in V}\frac{b(u,v)}{\|u\|\|v\|}=\gamma > 0$ can be defined in abstract infinite-dimensional mathematical spaces $U$ (where we seek the solution) and $V$ (which we use for testing in the Galerkin method),
 or in finite-dimensional spaces $U_h$ (where we seek the approximate solution) and $V_h$ (from where we select the test functions to generate a finite-dimensional system of linear equations to solve), namely
 $\inf_{u_h\in  U_h} \sup_{v_h \in V_h}\frac{b(u_h,v_h)}{\|u_h\|\|v_h\|}=\gamma_h > 0$. It may happen that even if we have the inf-sup  condition satisfied over the infinite-dimensional spaces $U,V$, when we will start to solve a problem in finite-dimensional subspaces $U_h,V_h$, the inf-sup condition will no longer be met. In this case 
 $\inf_{u_h\in  U_h} \sup_{v_h \in V_h}\frac{b(u_h,v_h)}{\|u_h\|\|v_h\|}=\gamma_h > \gamma = \inf_{u\in  U} \sup_{v \in V}\frac{b(u,v)}{\|u\|\|v\|}$. 
The problem with the inf-sup condition is that by selecting a finite-dimensional subspace of test functions $v_h \in V_h$, the supremum may no longer be controlled for all $u_h\in U_h$, by the same bound as in the infinite-dimensional case. For a large class of problems solved on computers using the finite element method, the continuous and discrete inf-sup condition is satisfied, and satisfactory solutions can be obtained. However, there is also a fairly large class of problems for which the discrete inf-sup condition is not met. The advection-dominated diffusion problem is an example of such a problem.
In such a case, the discrete solution obtained will not always be correct, e.g. numerical oscillations may arise. In that situation, it will be necessary to stabilize our problem. Some  strategies to obtain stable solutions in those cases are:
\begin{itemize}
\item Adding stabilizing terms to the  bilinear form
$b(u_h,v_h)$ such that those terms are zero in the ideal continuous finite-dimensional case, but they are non-zero in the discrete case. The idea is to improve the stability of the problem. An example of this method is the stabilization of the Streamline-Upwind Petrov-Galerkin method (SUPG) \cite{c4};
\item Modifying the equations of the problem and the discrete space where the solution is sought, in such a way that the inf-sup condition is obtained. This discrete space may not be conforming with respect to the infinite-dimensional space. An example of this strategy is the Discontinuous Galerkin (DG) method \cite{c5}, where functions defined on a mesh are discontinuous, and the equations are modified by adding terms that are equal to zero in infinite dimensional space, but not equal to zero in the discrete space. Then, a discrete norm is chosen so the inf-sup condition is satisfied at the  finite-dimensional level;
\item Reformulating the problem such that the approximation space and the test space are not the same, and imposing that the discrete test space is rich enough to obtain a stable method (towards the infinite-dimensional test space where the supremum is met). An example of this method is the residual minimization method \cite{c6}.
\end{itemize}
In this paper, we are interested in solving Problem~\eqref{eq:problem1} using a space-time formulation. Unlike time-stepping techniques, instead of generating a sequence of problems at each time step, we will solve on a computational space-time domain, where $x$ and $y$ are axes corresponding to the spatial dimensions and $t$ is the axis that corresponds to the temporal dimension. Thus, we will have a space-time discrete solution over a space-time mesh. The numerical solution for the space-time advection-dominated diffusion problem will also depend on the space-time method chosen. 

There are many recent research papers studying the feasibility of space-time formulations, i.e. considering time as another coordinate of the domain. In \cite{st3}, a summary of fast solvers for space-time formulations and different applications of the space-time methods can be found. More related to this work, a solver using space-time discretizations for the constrained first-order system least-square method (CFOSLS) is presented in \cite{st2}. Some space-time adaptive strategies can be found in  \cite{st1} and \cite{st4}. In \cite{st1}, numerical results for a space-time adaptive constrained first-order system with the least squares method are presented.
Paper \cite{st4} presents a comparison of algebraic multigrid methods for an adaptive space-time finite-element discretization of the heat equation in 3D and 4D. A space-time Discontinuous Petrov-Galerkin method is presented in~\cite{st5} for the second order Schrödinger equation.
Papers \cite{st6,st7} discuss the feasibility of the space-time Discontinuous Petrov-Galerkin formulation for the problem of acoustic wave propagation in multiple dimensions and show space-time adaptive mesh refinement strategies derived from the formulation.
Among others, we highlight the paper in \cite{st8} which presents the theory of a stable least-square finite element formulation for a parabolic space-time problem.

We develop a space-time solver to obtain numerical results for two least-square formulations on finite element meshes. The first formulation is based on isogeometric analysis (IGA) over  tensor product grids, which allows for fast (linear cost) inversion of the mass matrix.
The second one is based on the adaptive finite element method using  tetrahedral meshes.
The cost of the solver is related to the following three factors
\begin{itemize}
    \item The cost of inversion of the space-time mass matrix, which for the IGA discretization has a linear cost, and for the adaptive FEM discretization can be done using an iterative solver.
    \item The cost of the iterative solver for the solution of the system of equations of the Uzawa type \cite{Uzawa}. 
\end{itemize}    
We illustrate our findings on an advection-dominated diffusion space-time model problem.

The remainder of this paper is organized as follows. In Section~\ref{Sec:formulation} we introduce the model problem and its first-order formulation. In Section ~\ref{Sec:numericalform}, we show some numerical results without stabilization. In Section~\ref{Sec:Stabilization}, we introduce a stabilization constraint in the formulation. In Section~\ref{Sec:Isogeometric} we derive the solver algorithm. In Section~\ref{Stabilization}, we provide numerical results with a stability constraint. Finally, in Section~\ref{Sec:conclusions} conclusions and future work are presented. 

\section{Space-time formulation of the advection-dominated diffusion problem}\label{Sec:formulation}

Consider $\Omega\subset\mathbb{R}^2$ to be an open bounded domain with Lipschitz boundary $\partial\Omega$ and $T>0$.  We are interested in solving the follwing space-time advection-diffusion problem: Find a function $\phi:\overline{\Omega}\times[0,T]\to \mathbb{R}$, such that 
\begin{align}
\label{eq:adv-diff-problem}
  \Dt \phi - \varepsilon\Delta \phi + \boldsymbol{\beta} \cdot \nabla \phi &= f, \quad \text{ in }\Omega\times(0,T), \\
  \phi(\mathbf{x}, 0) &= u_0, \quad \text{in } \Omega\times\{0\}, \nonumber\\
  \phi(\mathbf{x}, t) &= 0, \quad \text{\,  in } \partial \Omega \times [0, T],\nonumber
\end{align}
for given functions $f:\Omega\times (0,T)\to \mathbb{R}$ and $u_0:\Omega\to \mathbb{R}$, and given parameters $\varepsilon\in \mathbb{R},$ and  $\boldsymbol{\beta}:=(\beta_x,\beta_y)^t\in \mathbb{R}^2$.  It is well-kown that for large~$\left\|\boldsymbol{\beta}\right\| / \varepsilon$, the problem is advection-dominated and the standard Galerkin method is unstable. 
A first-order system of equation equivalent to (\ref{eq:adv-diff-problem}) can be obtained by introducing a new variable $\SigU:=(\Lop\phi, \phi)^t$, where
\begin{equation}
\Lop\phi:=-\varepsilon\nabla\phi+\boldsymbol{\beta}\phi. \label{eq:Loperator}
\end{equation}
From the definition of $\SigU$ and ~\eqref{eq:adv-diff-problem} we have that $ \DivF \SigU = f.$

For sake of notation, we define $\Omega_T:=\Omega\times (0,T)$ to be the space-time domain and set $\Lop_x\phi:=-\varepsilon\partial_x\phi +\beta_x \phi$, and $\Lop_y\phi:=-\varepsilon\partial_y\phi +\beta_y \phi$. Thus $\Lop\phi:= (\Lop_x\phi, \Lop_y\phi)^t$. These notations and the selection of the auxiliary variable are motivated from the paper in~\cite{st2}. Thus, the problem presented in  ~\eqref{eq:adv-diff-problem} can be written as the first-order system: Find~$\SigU \in H(\DivF, \Omega_T)$, $\phi \in V$ such that
\begin{equation}
\left\{
\begin{aligned}
\SigU - \begin{bmatrix}\Lop\phi \\ \phi\end{bmatrix} &= 0, \\
\DivF \SigU &= f,
\end{aligned}
\right. \label{eq:first-order}
\end{equation}
where $V := \{ v \in H^1(\Omega_T) : \left.v\right|_{\Gamma_S} = 0, \,  v(\mathbf{x},0)= u_0(\mathbf{x}) \, a.e.\text{ in } \Omega\}$. Moreover, a simplified version of ~\eqref{eq:first-order} we can obtained, after imposing strongly ~$\SigU := (\Sig, \phi)$, i.e. 
\begin{equation}
\left\{
\begin{aligned}
\Dt \phi + \DivX \Sig &= f, \\
\Sig - \Lop\phi &= 0. \\
\end{aligned}
\right. \label{eq:reduced}
\end{equation}
A strong first-order variational formulation of  can be obtained by multiplying the first equation of~\eqref{eq:reduced} with a test function~$\psi \in L^2(\Omega_T)$ and integrating over $\Omega_T$, 
\begin{equation}
  \Prod{\Dt\phi}{\psi}_{L^2} + \Prod{\DivX\Sig}{\psi}_{L^2} = \Prod{f}{\psi}_{L^2}
\end{equation}
and multiplying the second equation of~\eqref{eq:reduced} with a test function~$\BTau := (\tau_x, \tau_y)^t \in L^2(\Omega_T) \times L^2(\Omega_T)$ and integrating over $\Omega$, to obtain 
\begin{equation}
  \Prod{\Sig}{\BTau}_{L^2} - \Prod{\Lop\phi}{\BTau}_{L^2} = 0.
\end{equation}
Thus, the formulation leads us to seek for~$\SigU \in H(\text{div}, \Omega_T)$ and $\phi \in V$ such that 
\begin{equation}
\left\{
\begin{aligned}
  \Prod{\Dt\phi}{\psi}_{L^2} + \Prod{\D{x}\sigma_x}{\psi}_{L^2} + \Prod{\D{y}\sigma_y}{\psi}_{L^2} &= \Prod{f}{\psi}_{L^2}, \\
  \Prod{\sigma_x}{\tau_x}_{L^2} - \Prod{\Lop_x \phi}{\tau_x}_{L^2} &= 0, \\
  \Prod{\sigma_y}{\tau_y}_{L^2} - \Prod{\Lop_y \phi}{\tau_y}_{L^2} &= 0. 
\end{aligned}
\right. \label{eq:strong}
\end{equation}
for all $\psi \in  V$ and $\boldsymbol{\tau}\in H(\text{div},\Omega_T)$. 

\vfill
\textbf{Discrete space-time formulation.} Consider a conforming mesh $\Omega_h$ of $\Omega$, and a let $I_h=\{I_1,I_2,...,I_N\}$ be a partition of the time interval $[0,T]$, with $I_k=[t_k, t_{k+1}]$, $h>0$,  and  $t_k= h(k-1),$ for $k=1,\dots N$. Then, define $\Omega_{T,h}:=\Omega_h\times I_h$ to be a tensor-product  mesh of $\Omega_T:=\Omega\times(0,T)$ and $\Gamma_S := \partial \Omega \times (0, T)$ to be the evolving-in-time spatial boundary. 

Selecting appropriated finite-dimensional spaces $\Sigma_h$ and $V_h$, subspaces of $ H(\text{div},\Omega_T)$ and $V$ respectively, the discrete first-order formulation derived from~\eqref{eq:strong} seeks a discrete solution $\phi^h\in \Sigma_h$ and $\Sig^h\in \Sigma_h$ such that  
\begin{equation}
\left\{
\begin{aligned}
  \Prod{\Dt\phi^h}{\psi}_{L^2} + \Prod{\D{x}\sigma_x^h}{\psi}_{L^2} + \Prod{\D{y}\sigma_y^h}{\psi^h}_{L^2} &= \Prod{f}{\psi^h}_{L^2}, \\
  \Prod{\sigma_x^h}{\tau_x^h}_{L^2} - \Prod{\Lop_x \phi^h}{\tau_x^h}_{L^2} &= 0, \\
  \Prod{\sigma_y^h}{\tau_y^h}_{L^2} - \Prod{\Lop_y \phi}{\tau_y^h}_{L^2} &= 0, \\
\end{aligned}
\right.
\end{equation}
for all $\psi\in V_h$ and $\boldsymbol{\tau}\in \Sigma_h$. In our case, we approximate the independent variables $\phi^h,\sigma^h_x,$ and $\sigma^h_y$ using quadratic B-splines. Let us define $\{u_i\}$ the B-spline basis functions and $\boldsymbol{\phi}^{\texttt{h}}_x,\boldsymbol{\sigma}^{\texttt{h}}_x, \boldsymbol{\sigma}^{\texttt{h}}_y,$  the corresponding vectors of coefficients of the $B$-spline expansion of $\phi^h, \sigma_x^h, \sigma_x^h$, respectively.
We can write the system in the following matrix structure
\begin{equation}
  \begin{bmatrix}
  \color{blue}{A_t} & \color{blue}{A_x} & \color{blue}{A_y} \\
  \color{cadmiumgreen}{L_x} & \color{purple}{M} & 0 \\
  \color{cadmiumgreen}{L_y} & 0 & \color{purple}{M} \\
  \end{bmatrix}
  \begin{bmatrix}
  \boldsymbol{\phi}^{\texttt{h}} \\ \boldsymbol{\sigma}_x^{\texttt{h}} \\ \boldsymbol{\sigma}_y^{\texttt{h}}
  \end{bmatrix}
  =
  \begin{bmatrix}
  \boldsymbol{f}^{\texttt{h}} \\ 0 \\ 0
  \end{bmatrix},
\end{equation}
where $\color{purple}{M}$ represents the mass matrix such that $\color{purple}{(M)_{ij}:=(u_i,u_j)_{L^2}}$,  $\color{blue}({A}_\gamma)$, is such that  $\color{blue}{(A_{\gamma})_{ij}=\Prod{\D{\gamma}u_i}{u_j}}$ and $\color{cadmiumgreen}{L_\gamma}$ such that $\color{cadmiumgreen}{(L_\gamma)_{ij}:=\Prod{\Lop_\gamma u_i}{u_j}}.$ Moreover, $\mathbf{f}^{\texttt{h}}$ is also the vector of coefficients related to the expansion of $f$ in the B-spline basis.

\section{Numerical results without stabilization}\label{Sec:numericalform}

\subsection{Space-time pure diffusion problem}

\begin{figure}[h]
\centering
  \begin{subfigure}[b]{0.32\textwidth}
  \centering
  \includegraphics[width=\textwidth]{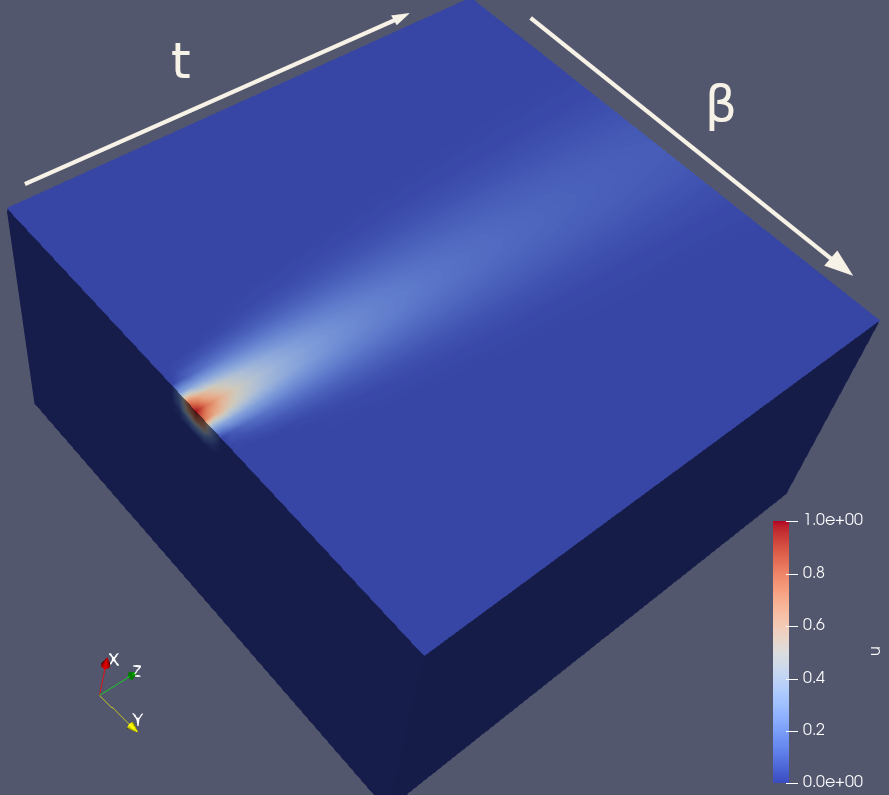}
  \caption{$\varepsilon = 5\times 10^{-3}$}
  \end{subfigure}
  \begin{subfigure}[b]{0.32\textwidth}
  \centering
  \includegraphics[width=\textwidth]{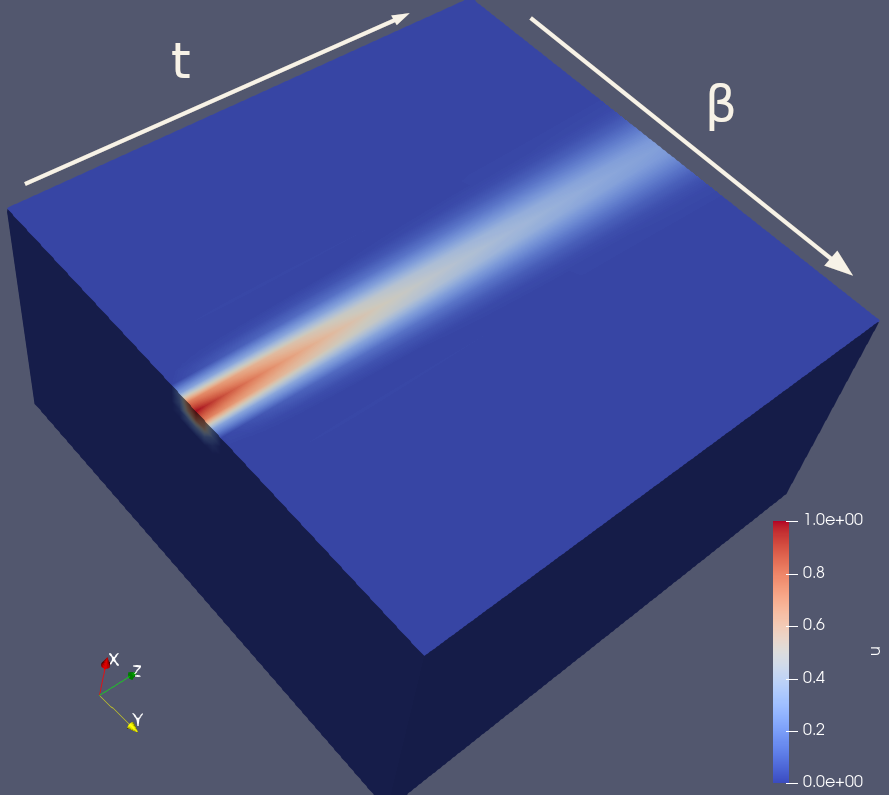}
  \caption{$\varepsilon = 10^{-3}$}
  \end{subfigure}
  \begin{subfigure}[b]{0.32\textwidth}
  \centering
  \includegraphics[width=\textwidth]{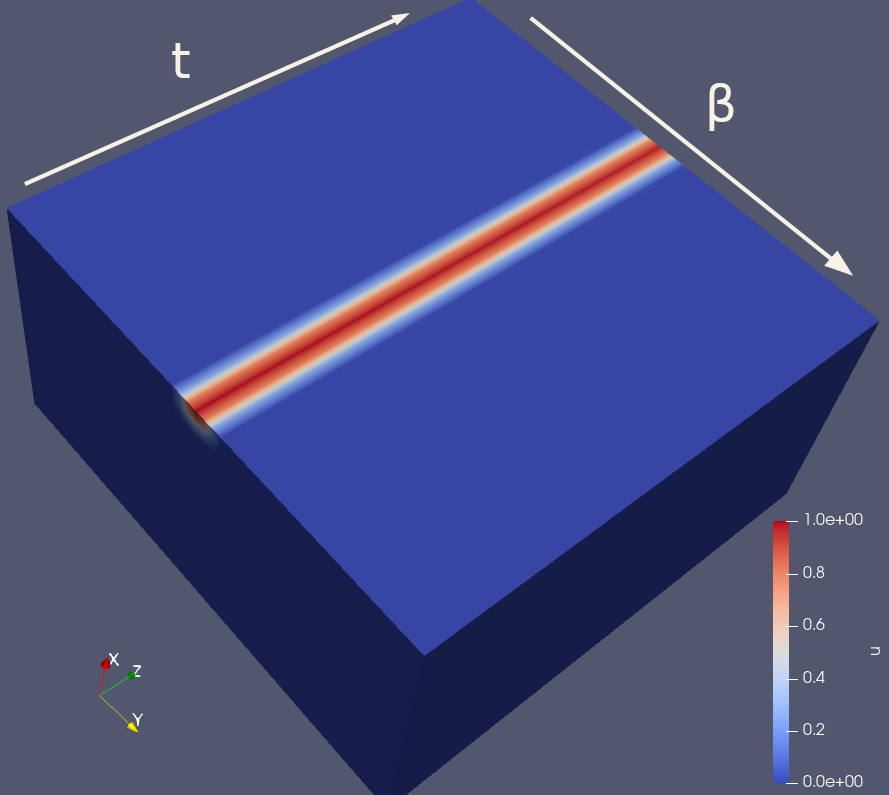}
  \caption{$\varepsilon = 10^{-5}$}
  \end{subfigure}

  \caption{Pure diffusion. Results for~$\boldsymbol{\beta} = \mathbf{0}$, $T = 1$, and  initial condition~\eqref{eq:init}, in a mesh with $32\times 32 \times 32$ elements, approximating the solution $\phi$ with quadratic B-splines.}
  \label{fig:purediff}
\end{figure}

The goal of the first numerical experiment is to verify the correctness of the implementation, and show that the problem without advection does not require stabilization.
The first numerical experiments concerns the model problem~\eqref{eq:adv-diff-problem}
formulated on a regular space-time domain~$\Omega_T= (0, 1)^3$
with the advection "wind" vector $\boldsymbol{\beta} = (0, 0)^t$, the diffusion coefficient
$\varepsilon = 10^{-5}$, and without a source force ($f = 0$). The initial state $u_0$ is given by~
 \begin{eqnarray}
 u_0(\mathbf{x}) = \psi\left(10\|\mathbf{x} - \mathbf{c}\|\right), \quad
  \psi(r) =
  \begin{cases}
     (1 - r^2)^2, & \text{ for } r \leq 1, \\
     0, & \text{ for } r > 1,
  \end{cases}
  \label{eq:init}
\end{eqnarray}
with $\mathbf{c} = (0.5, 0.5)^t$. 
As a result, the initial state is zero except for a small region in the center of the domain.
Thus, in order to invert the space-time stiffness matrix during our solver algorithm, we assume that the solution is indeed zero at some point, e.g. $(0,0,0)^t$.
We can read from Figure \ref{fig:purediff} that our code provides a stable numerical solution.

\subsection{Space-time advection-diffusion problem}

\begin{figure}[h]
\centering
  \begin{subfigure}[b]{0.32\textwidth}
  \centering
  \includegraphics[width=\textwidth]{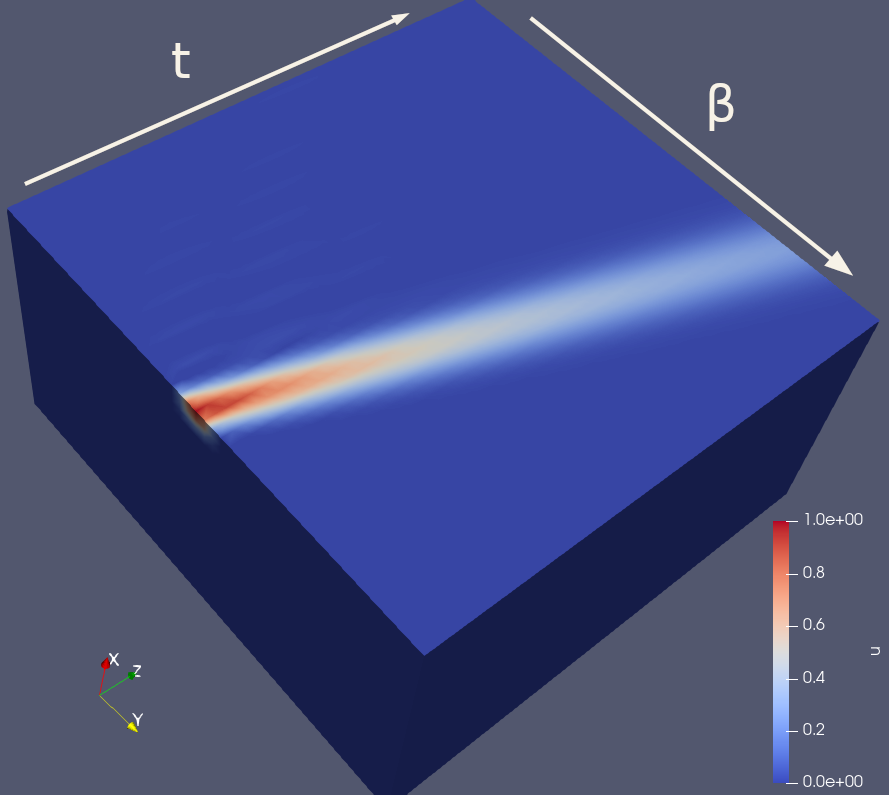}
  \caption{$\varepsilon = 10^{-3}$, $s = 0.3$}
  \end{subfigure}
  \begin{subfigure}[b]{0.32\textwidth}
  \centering
  \includegraphics[width=\textwidth]{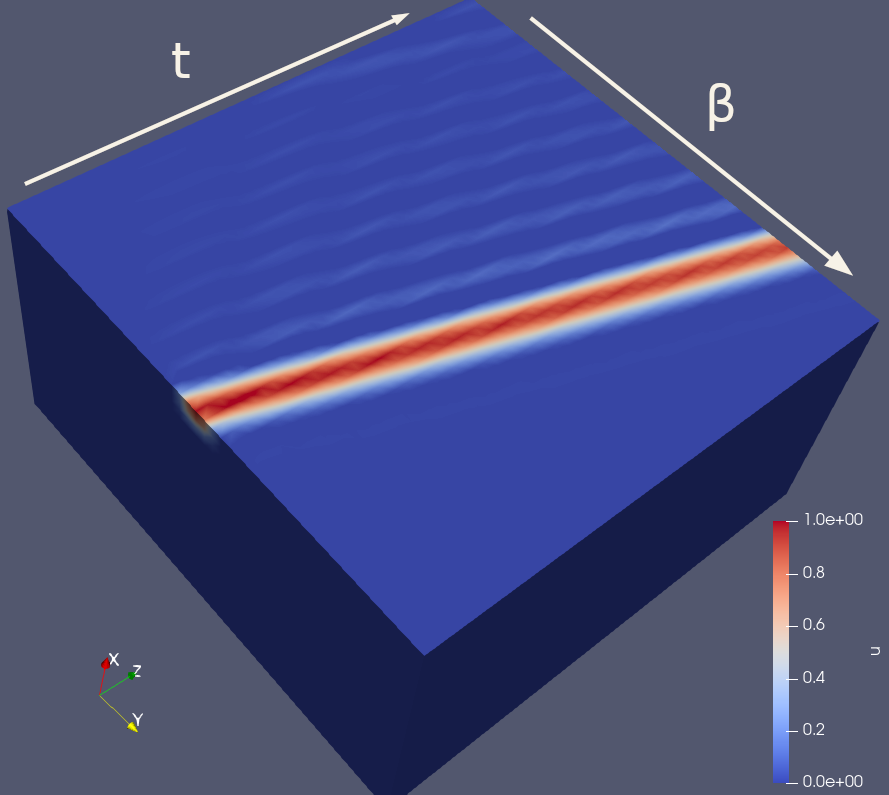}
  \caption{$\varepsilon = 10^{-5}$, $s = 0.3$}
  \end{subfigure}
  \begin{subfigure}[b]{0.32\textwidth}
  \centering
  \includegraphics[width=\textwidth]{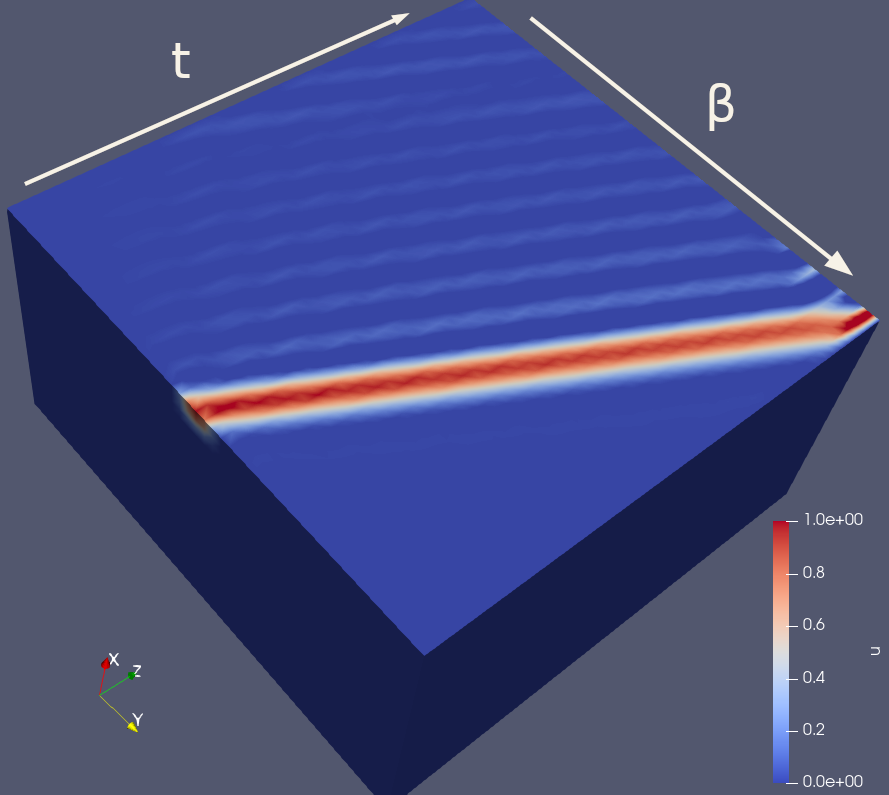}
  \caption{$\varepsilon = 10^{-6}$, $s = 0.5$}
  \end{subfigure}

  \begin{subfigure}[b]{0.32\textwidth}
  \includegraphics[width=\textwidth]{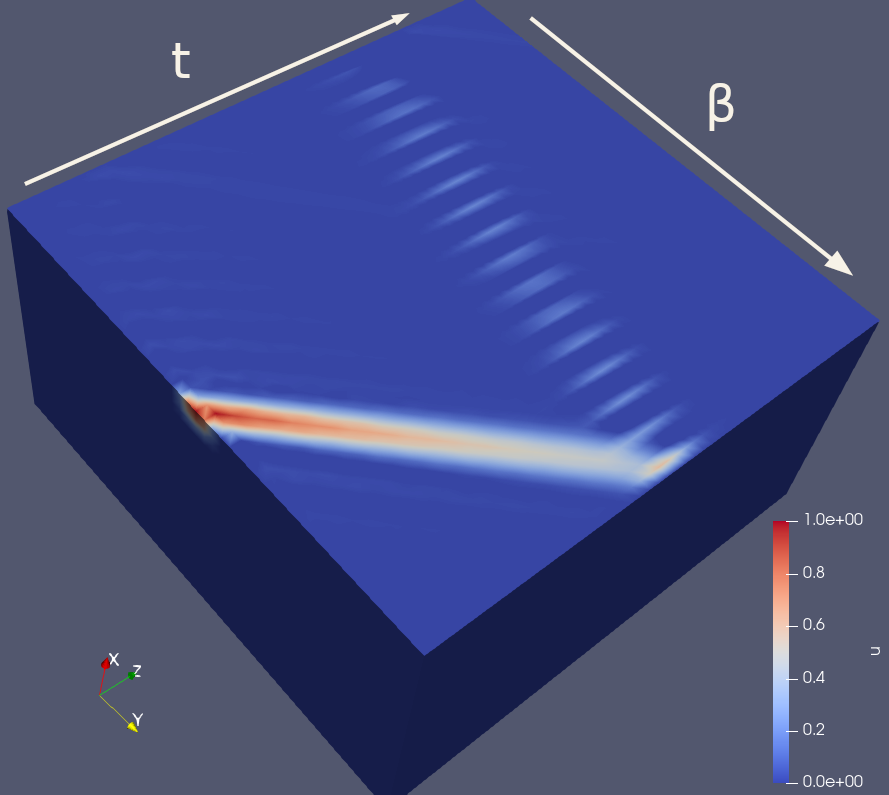}
  \caption{$\varepsilon = 10^{-3}$, $s = 1$}
  \end{subfigure}
  \begin{subfigure}[b]{0.32\textwidth}
  \centering
  \includegraphics[width=\textwidth]{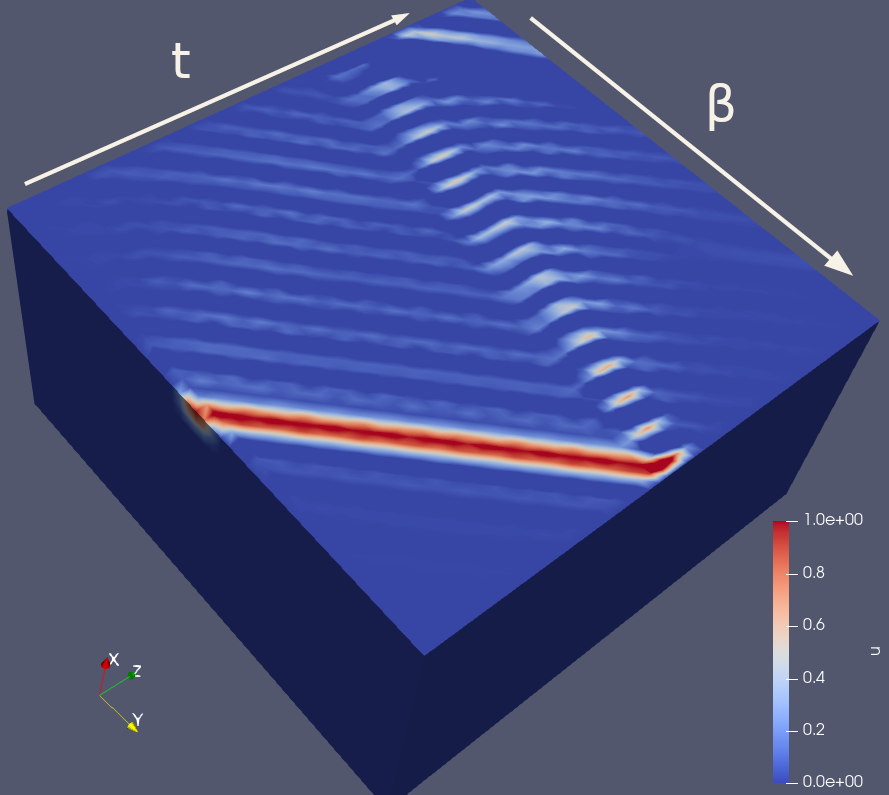}
  \caption{$\varepsilon = 10^{-5}$, $s = 1$}
  \end{subfigure}
  \begin{subfigure}[b]{0.32\textwidth}
  \centering
  \includegraphics[width=\textwidth]{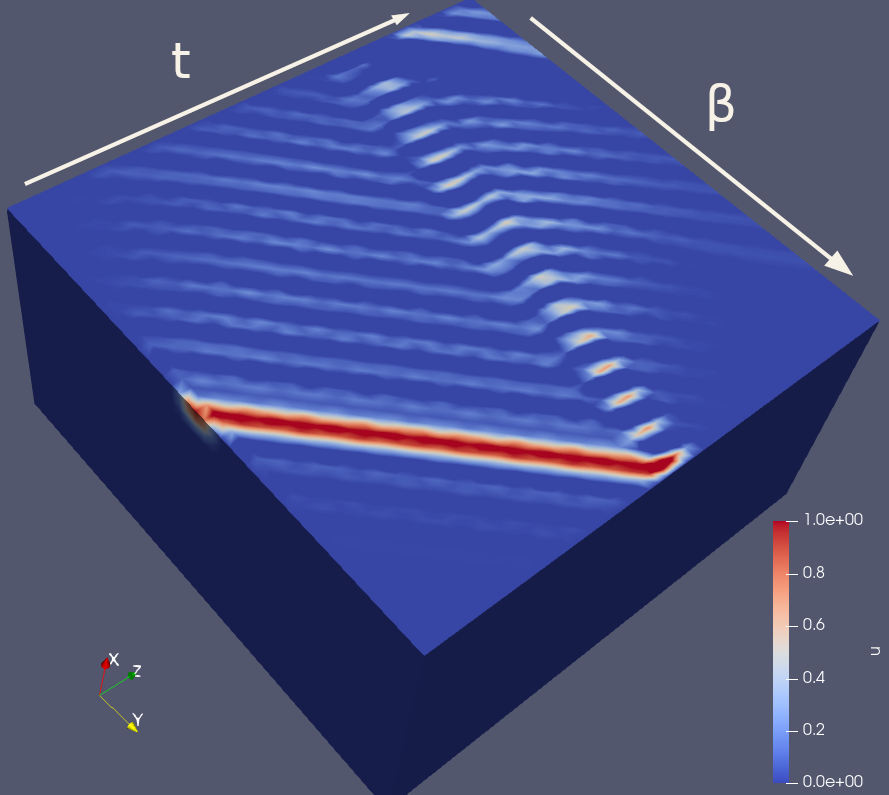}
  \caption{$\varepsilon = 10^{-6}$, $s = 1$}
  \end{subfigure}

  \caption{Advection-diffusion. Results for~$\beta = (0, s)$, $T = 1$,  and  initial condition~\eqref{eq:init} in a mesh with $32\times 32 \times 32$ elements, using quadratic B-splines.}
  \label{fig:advdiff}
\end{figure}

The second numerical experiment for the model problem~\eqref{eq:adv-diff-problem}
formulated on a regular space-time domain~$\Omega_T= (0, 1)^3$ introduces the advection "wind" vector $\boldsymbol{\beta} = (0, 0.3)^t$.
The diffusion coefficient
$\varepsilon = 10^{-5}$, and 
 forcing $f = 0$. The initial state $u_0$ is given by (\ref{eq:init}).
The initial state is zero except for a small region in the center of the domain. The numerical results summarized in Figure \ref{fig:advdiff} are unstable. They introduce unexpected oscillations that grow when $\left\|\boldsymbol{\beta}\right\| / \varepsilon$ gets bigger.
We need to introduce our stabilization method.

\section{Constraint Least-Square Stabilized method}\label{Sec:Stabilization}

In this section, we introduce a stabilization method for the space-time advection-dominated diffusion problem based on the  constrained least-squares method, arising from the first-order constraint minimization problem
\begin{equation}
\begin{aligned}
  \textbf{min} \ & J(\SigU, \phi), \\
  \textbf{subject to} \  & \DivF \SigU = f,
\end{aligned}
\label{eq:stab}
\end{equation}
where $J:H(\DivF, \Omega_T)\times V\to \mathbb{R}^+_0$ is defined by $$J(\SigU,\phi):=\frac{1}{2} \left\| \SigU - \begin{bmatrix} \Lop\phi \\ \phi \end{bmatrix}\right\|^2_{L^2}
  =
  \frac{1}{2} \left\| \Sig - \Lop\phi \right\|^2_{L^2}
  + \frac{1}{2} \left\| \sigma_* - \phi\right\|^2_{L^2},$$ $\Lop\phi$ is defined as in~\eqref{eq:Loperator} and~$\SigU := (\Sig, \sigma_*)^t$.
Introducing the Lagrange multiplier function $\lambda$, this problem is equivalent to minimize the functional $G$, defined by
\begin{equation}
  G(\SigU, \phi, \lambda) =
  \frac{1}{2} \left\| \Sig - \Lop\phi \right\|^2_{L^2} +
  \frac{1}{2} \left\| \sigma_* - \phi\right\|^2_{L^2} +
  {\Prod{\DivF \SigU - f}{\lambda}}_{L^2}.
\end{equation}
%
The equivalent problem reads as follows:
find~$(\SigU$, $\phi$, $\lambda) \in \mathbf{W}$
such that for all~$(\BTauU$, $\omega$, $\mu) \in \mathbf{W}$
\begin{equation}
\begin{aligned}
  \Prod{\Sig - \Lop\phi}{\BTau - \Lop\omega}_{L^2} + \Prod{\sigma_* - \phi}{\tau_* - \omega}_{L^2} \\+
  \Prod{\DivF\BTauU}{\lambda}_{L^2} + \Prod{\DivF \SigU - f}{\mu}_{L^2} = 0,
\end{aligned} \label{eq:system}
\end{equation}
where~$\mathbf{W} := H(\text{div}, \Omega_T) \times V \times L^2(\Omega_T)$.

The resulting saddle-point system derived from \eqref{eq:system} reads: Find $\phi,\SigU,\lambda\in \mathbf{W}$ such that 
\begin{equation}
\begin{alignedat}{6}
  &\Prod{\phi}{\omega}_{L^2} + \Prod{\Lop\phi}{\Lop\omega}_{L^2} &&- \Prod{\sigma_*}{\omega}_{L^2} &&- \Prod{\Sig}{\Lop\omega}_{L^2} && &&=0,\\
  &\hspace{5mm}-\Prod{\phi}{\tau_*}_{L^2} &&+ \Prod{\sigma_*}{\tau_*}_{L^2} && &&+ \Prod{\lambda}{\D{t}\tau_*}_{L^2} &&=0, \\
  &\hspace{5mm}-\Prod{\Lop\phi}{\BTau}_{L^2} && &&+ \Prod{\Sig}{\BTau}_{L^2} &&+ \Prod{\lambda}{\DivX\BTau}_{L^2} &&=0, \\
  & &&+\Prod{\D{t}\sigma_*}{\mu}_{L^2} &&+ \Prod{\DivX\Sig}{\mu}_{L^2} && &&=\Prod{f}{\mu}_{L^2},
\end{alignedat} \label{eq:constraint-problem}
\end{equation}
for all $\omega, \BTauU, \mu\in \mathbf{W}.$

In what follows, we study the feasibility of a solver for the IGA space-time discretization on tensor product grids for the space-time formulation in~\eqref{eq:constraint-problem}. The discrete problem arising from the system in ~\eqref{eq:constraint-problem} can be expressed in its matrix form as follows:
\begin{equation}
  \begin{bmatrix}
  \color{purple}{M_{\phi}} + \color{red}{K_{\phi}} & -\color{purple}{M^T_{\phi\sigma}} & -\color{cadmiumgreen}{L_x^T} & -\color{cadmiumgreen}{L_y^T} & 0 \\
  -\color{purple}{M_{\phi\sigma}} & \color{purple}{M} & 0 & 0 & \color{blue}{A_t^T} \\
  -\color{cadmiumgreen}{L_x} & 0 & \color{purple}{M} & 0 & \color{blue}{A_x^T} \\
  -\color{cadmiumgreen}{L_y} & 0 & 0 & \color{purple}{M} & \color{blue}{A_y^T} \\
  0 & \color{blue}{A_t} & \color{blue}{A_x} & \color{blue}{A_y} & 0 \\
  \end{bmatrix}
  \begin{bmatrix}
  \phi^{h} \\ \sigma_*^{h} \\ \sigma_x^{h} \\ \sigma_y^{h} \\ \lambda^{h}
  \end{bmatrix}
  =
  \begin{bmatrix}
  0 \\ 0 \\ 0 \\ 0 \\ f^{h}
  \end{bmatrix},
\end{equation}
where
\begin{align}
\color{purple}{M}_{ijk,lmn}&:=(B^x_iB^y_jB^t_k,B^x_lB^y_mB^t_n)_{L^2}, \notag \\
\color{blue}{A_x}_{ijk,lmn}&:=((\partial_xB^x_i)B^y_jB^t_k,B^x_lB^y_mB^t_n)_{L^2}, \notag \\
\color{blue}{A_y}_{ijk,lmn}&:=(B^x_i(\partial_yB^y_j)B^t_k,B^x_lB^y_mB^t_n)_{L^2}, \notag \\
\color{blue}{A_t}_{ijk,lmn}&:=(B^x_iB^y_j(\partial_tB^t_k),B^x_lB^y_mB^t_n)_{L^2}, \notag \\
\color{cadmiumgreen}{L_x}_{ijk,lmn}&:=(({\cal L}_xB^x_i)B^y_jB^t_k,B^x_lB^y_mB^t_n)_{L^2}, \notag \\
\color{cadmiumgreen}{L_y}_{ijk,lmn}&:=(B^x_i({\cal L}_yB^y_j)B^t_k,B^x_lB^y_mB^t_n)_{L^2}, \notag \\
\color{red}{K}_{ijk,lmn}&:=(({\cal L}_xB^x_i)B^y_jB^t_k,({\cal L}_xB^x_l)B^y_mB^t_n)_{L^2}, +
(B^x_i({\cal L}_yB^y_j)B^t_k,B^x_l({\cal L}_yB^y_m)B^t_n)_{L^2}, \notag \\
{S}_{ijk,lmn}&:=\color{blue}{A_x}_{ijk,lmn}+\color{blue}{A_y}_{ijk,lmn}+
\color{blue}{A_t}_{ijk,lmn}. 
\label{eq:matrices}
\end{align}
\revision{Here, the ${\color{purple}{M_{\phi}}}$ denotes the mass matrix including the Dirichlet boundary conditions for the concentration field, and ${\color{purple}{M_{\phi\sigma}}}$ denotes the mass matrix with the trial space including the zero Dirichlet b.c., and the test space without the zero Dirichlet b.c.}

\section{Solver algorithm for the space-time formulation}\label{Sec:Isogeometric}
We permute the matrix and obtained the equivalent matrix problem:
\begin{equation}
  \begin{bmatrix}
  \color{purple}{M} & 0 & 0 & \color{blue}{A_t^T} & -\color{purple}{M^T_{\phi\sigma}}  \\
   0 & \color{purple}{M} & 0 & \color{blue}{A_x^T} & -\color{cadmiumgreen}{L_x}  \\
  0 & 0 & \color{purple}{M} & \color{blue}{A_y^T} & -\color{cadmiumgreen}{L_y} \\
  \color{blue}{A_t} & \color{blue}{A_x} & \color{blue}{A_y} & 0 & 0 \\
  -\color{purple}{M_{\phi\sigma}} & -\color{cadmiumgreen}{L_x^T} & -\color{cadmiumgreen}{L_y^T} & 0 & \color{purple}{M_{\phi}} + \color{red}{K} \\
  \end{bmatrix}
  \begin{bmatrix}
  \sigma_* \\ \sigma_x \\ \sigma_y \\ \lambda \\ \phi 
  \end{bmatrix}
  =
  \begin{bmatrix}
  0 \\ 0 \\ 0 \\ f \\ 0
  \end{bmatrix}.
\end{equation}
We solve the system in the following steps:
\begin{enumerate}
\item We compute $\color{red}{M}^{-1}$. For IGA discretization on tensor product grids this leads to $\color{red}{M}^{-1}=\color{red}{M_x}^{-1}\otimes \color{red}{M_y}^{-1} \otimes \color{red}{M_t}^{-1}$ which can be done in a linear ${\cal O}(N)$ cost due to the Kronecker product structure of $\color{red}{M}$.
\item We multiple the first row by $\color{purple}{M}^{-1}$
\begin{equation}
  \begin{bmatrix}
  {I} & 0 & 0 & \color{purple}{M}^{-1}\color{blue}{A_t^T} & -\color{purple}{M^{-1}M^T_{\phi\sigma}}  \\
   0 & \color{purple}{M} & 0 & \color{blue}{A_x^T} & -\color{cadmiumgreen}{L_x}  \\
  0 & 0 & \color{purple}{M} & \color{blue}{A_y^T} & -\color{cadmiumgreen}{L_y} \\
  \color{blue}{A_t} & \color{blue}{A_x} & \color{blue}{A_y} & 0 & 0 \\
  -\color{purple}{M_{\phi\sigma}} & -\color{cadmiumgreen}{L_x^T} & -\color{cadmiumgreen}{L_y^T} & 0 & \color{purple}{M_{\phi}} + \color{red}{K} \\
  \end{bmatrix}
  \begin{bmatrix}
  \sigma_* \\ \sigma_x \\ \sigma_y \\ \lambda \\ \phi 
  \end{bmatrix}
  =
  \begin{bmatrix}
  0 \\ 0 \\ 0 \\ f \\ 0
  \end{bmatrix}
\end{equation}\item We subtract the first row multiplied by $\color{blue}{A_t}$ from the fourth row
\begin{equation}
  \begin{bmatrix}
  {I} & 0 & 0 & \color{purple}{M}^{-1}\color{blue}{A_t^T} & -\color{purple}{M^{-1}M^T_{\phi\sigma}}  \\
   0 & \color{purple}{M} & 0 & \color{blue}{A_x^T} & -\color{cadmiumgreen}{L_x}  \\
  0 & 0 & \color{purple}{M} & \color{blue}{A_y^T} & -\color{cadmiumgreen}{L_y} \\
  0 & \color{blue}{A_x} & \color{blue}{A_y} & -\color{blue}{A_t}\color{purple}{M}^{-1}\color{blue}{A_t^T} & \color{blue}{A_t}\color{purple}{M^{-1}M^T_{\phi\sigma}} \\
  -\color{purple}{M_{\phi\sigma}} & -\color{cadmiumgreen}{L_x^T} & -\color{cadmiumgreen}{L_y^T} & 0 & \color{purple}{M_{\phi}} + \color{red}{K} \\
  \end{bmatrix}
  \begin{bmatrix}
  \sigma_* \\ \sigma_x \\ \sigma_y \\ \lambda \\ \phi 
  \end{bmatrix}
  =
  \begin{bmatrix}
  0 \\ 0 \\ 0 \\ f \\ 0
  \end{bmatrix} 
  \end{equation}
  \item We add the first row multiplied by $\color{purple}{M_{\phi\sigma}}$ to the fifth row
\begin{equation}
  \begin{bmatrix}
  {I} & 0 & 0 & \color{purple}{M}^{-1}\color{blue}{A_t^T} & -\color{purple}{M^{-1}M^T_{\phi\sigma}}  \\
   0 & \color{purple}{M} & 0 & \color{blue}{A_x^T} & -\color{cadmiumgreen}{L_x}  \\
  0 & 0 & \color{purple}{M} & \color{blue}{A_y^T} & -\color{cadmiumgreen}{L_y} \\
  0 & \color{blue}{A_x} & \color{blue}{A_y} & -\color{blue}{A_t}\color{purple}{M}^{-1}\color{blue}{A_t^T} & \color{blue}{A_t}\color{purple}{M^{-1}M^T_{\phi\sigma}} \\
  0 & -\color{cadmiumgreen}{L_x^T} & -\color{cadmiumgreen}{L_y^T} & \color{purple}{M_{\phi\sigma}}\color{purple}{M}^{-1}\color{blue}{A_t^T} & 
  -\color{purple}{M_{\phi\sigma} M^{-1}M^T_{\phi\sigma}}+\color{purple}{M_{\phi}} + \color{red}{K} \\
  \end{bmatrix}
  \begin{bmatrix}
  \sigma_* \\ \sigma_x \\ \sigma_y \\ \lambda \\ \phi 
  \end{bmatrix}
  =
  \begin{bmatrix}
  0 \\ 0 \\ 0 \\ f \\ 0
  \end{bmatrix} 
  \end{equation}
\item We multiple the second row by $\color{purple}{M}^{-1}$
\begin{equation}
  \begin{bmatrix}
  {I} & 0 & 0 & \color{purple}{M}^{-1}\color{blue}{A_t^T} & -\color{purple}{M^{-1}M^T_{\phi\sigma}}  \\
   0 & I & 0 & \color{purple}{M^{-1}}\color{blue}{A_x^T} & -\color{purple}{M^{-1}}\color{cadmiumgreen}{L_x}  \\
  0 & 0 & \color{purple}{M} & \color{blue}{A_y^T} & -\color{cadmiumgreen}{L_y} \\
  0 & \color{blue}{A_x} & \color{blue}{A_y} & -\color{blue}{A_t}\color{purple}{M}^{-1}\color{blue}{A_t^T} & \color{blue}{A_t}\color{purple}{M^{-1}M^T_{\phi\sigma}} \\
  0 & -\color{cadmiumgreen}{L_x^T} & -\color{cadmiumgreen}{L_y^T} & \color{purple}{M_{\phi\sigma}}\color{purple}{M}^{-1}\color{blue}{A_t^T} & 
  -\color{purple}{M_{\phi\sigma} M^{-1}M^T_{\phi\sigma}}+\color{purple}{M_{\phi}} + \color{red}{K} \\
  \end{bmatrix}
  \begin{bmatrix}
  \sigma_* \\ \sigma_x \\ \sigma_y \\ \lambda \\ \phi 
  \end{bmatrix}
  =
  \begin{bmatrix}
  0 \\ 0 \\ 0 \\ f \\ 0
  \end{bmatrix} 
  \end{equation}
    \item We subtract the second row multiplied by $\color{blue}{A_x}$ from the fourth row
\begin{eqnarray}
  \begin{bmatrix}
  {I} & 0 & 0 & \color{purple}{M}^{-1}\color{blue}{A_t^T} & -\color{purple}{M^{-1}M^T_{\phi\sigma}}  \\
   0 & I & 0 & \color{purple}{M^{-1}}\color{blue}{A_x^T} & -\color{purple}{M^{-1}}\color{cadmiumgreen}{L_x}  \\
  0 & 0 & \color{purple}{M} & \color{blue}{A_y^T} & -\color{cadmiumgreen}{L_y} \\
  0 & 0 & \color{blue}{A_y} & -(\color{blue}{A_t}\color{purple}{M}^{-1}\color{blue}{A_t^T}+
  \color{blue}{A_x}\color{purple}{M}^{-1}\color{blue}{A_x^T}) & (\color{blue}{A_t}\color{purple}{M^{-1}M^T_{\phi\sigma}}+\color{blue}{A_x}\color{purple}{M}^{-1}\color{cadmiumgreen}{L_x}) \\
  0 & -\color{cadmiumgreen}{L_x^T} & -\color{cadmiumgreen}{L_y^T} & \color{purple}{M_{\phi\sigma}}\color{purple}{M}^{-1}\color{blue}{A_t^T} & 
  -\color{purple}{M_{\phi\sigma} M^{-1}M^T_{\phi\sigma}}+\color{purple}{M_{\phi}} + \color{red}{K} \\  \end{bmatrix} \notag \\
  \begin{bmatrix}
  \sigma_* \\ \sigma_x \\ \sigma_y \\ \lambda \\ \phi 
  \end{bmatrix}
  =
  \begin{bmatrix}
  0 \\ 0 \\ 0 \\ f \\ 0
  \end{bmatrix} 
  \end{eqnarray}
      \item We add the second row multiplied by $\color{cadmiumgreen}{L^T_x}$ to the fifth row
\begin{eqnarray}
  \begin{bmatrix}
  {I} & 0 & 0 & \color{purple}{M}^{-1}\color{blue}{A_t^T} & -\color{purple}{M^{-1}M^T_{\phi\sigma}}  \\
   0 & I & 0 & \color{purple}{M^{-1}}\color{blue}{A_x^T} & -\color{purple}{M^{-1}}\color{cadmiumgreen}{L_x}  \\
  0 & 0 & \color{purple}{M} & \color{blue}{A_y^T} & -\color{cadmiumgreen}{L_y} \\
  0 & 0 & \color{blue}{A_y} & -(\color{blue}{A_t}\color{purple}{M}^{-1}\color{blue}{A_t^T}+
  \color{blue}{A_x}\color{purple}{M}^{-1}\color{blue}{A_x^T}) & (\color{blue}{A_t}\color{purple}{M^{-1}M^T_{\phi\sigma}}+\color{blue}{A_x}\color{purple}{M}^{-1}\color{cadmiumgreen}{L_x}) \\
  0 & 0 & -\color{cadmiumgreen}{L_y^T} & \color{cadmiumgreen}{L_x^T} \color{purple}{M^{-1}}\color{blue}{A_x^T}
 +\color{purple}{M_{\phi\sigma}}\color{purple}{M}^{-1}\color{blue}{A_t^T} & 
 -\color{cadmiumgreen}{L_x^T}\color{purple}{M^{-1}}\color{cadmiumgreen}{L_x}
  -\color{purple}{M_{\phi\sigma} M^{-1}M^T_{\phi\sigma}}+\color{purple}{M_{\phi}} + \color{red}{K} \\
  \end{bmatrix} \notag \\
    \begin{bmatrix}
  \sigma_* \\ \sigma_x \\ \sigma_y \\ \lambda \\ \phi 
  \end{bmatrix}
  =
  \begin{bmatrix}
  0 \\ 0 \\ 0 \\ f \\ 0
  \end{bmatrix} \notag \\
  \end{eqnarray}
\item We multiple the third row by $\color{purple}{M}^{-1}$

\begin{eqnarray}
  \begin{bmatrix}
  {I} & 0 & 0 & \color{purple}{M}^{-1}\color{blue}{A_t^T} & -\color{purple}{M^{-1}M^T_{\phi\sigma}}  \\
   0 & I & 0 & \color{purple}{M^{-1}}\color{blue}{A_x^T} & -\color{purple}{M^{-1}}\color{cadmiumgreen}{L_x}  \\
  0 & 0 & I & \color{purple}{M^{-1}}\color{blue}{A_y^T} & -\color{purple}{M^{-1}}\color{cadmiumgreen}{L_y} \\
  0 & 0 & \color{blue}{A_y} & -(\color{blue}{A_t}\color{purple}{M}^{-1}\color{blue}{A_t^T}+
  \color{blue}{A_x}\color{purple}{M}^{-1}\color{blue}{A_x^T}) & (\color{blue}{A_t}\color{purple}{M^{-1}M^T_{\phi\sigma}}+\color{blue}{A_x}\color{purple}{M}^{-1}\color{cadmiumgreen}{L_x}) \\
  0 & 0 & -\color{cadmiumgreen}{L_y^T} & \color{cadmiumgreen}{L_x^T} \color{purple}{M^{-1}}\color{blue}{A_x^T}
 +\color{purple}{M_{\phi\sigma}}\color{purple}{M}^{-1}\color{blue}{A_t^T} & 
 -\color{cadmiumgreen}{L_x^T}\color{purple}{M^{-1}}\color{cadmiumgreen}{L_x}
  -\color{purple}{M_{\phi\sigma} M^{-1}M^T_{\phi\sigma}}+\color{purple}{M_{\phi}} + \color{red}{K} \\
  \end{bmatrix} \notag \\
    \begin{bmatrix}
  \sigma_* \\ \sigma_x \\ \sigma_y \\ \lambda \\ \phi 
  \end{bmatrix}
  =
  \begin{bmatrix}
  0 \\ 0 \\ 0 \\ f \\ 0
  \end{bmatrix} \notag \\
  \end{eqnarray}
      \item We subtract the third row multiplied by $\color{blue}{A_y}$ from the fourth row
{\footnotesize{
\begin{eqnarray}
  \begin{bmatrix}
  {I} & 0 & 0 & \color{purple}{M}^{-1}\color{blue}{A_t^T} & -\color{purple}{M^{-1}M^T_{\phi\sigma}}  \\
   0 & I & 0 & \color{purple}{M^{-1}}\color{blue}{A_x^T} & -\color{purple}{M^{-1}}\color{cadmiumgreen}{L_x}  \\
  0 & 0 & I & \color{purple}{M^{-1}}\color{blue}{A_y^T} & -\color{purple}{M^{-1}}\color{cadmiumgreen}{L_y} \\
  0 & 0 & 0 & 
  -(\color{blue}{A_y}\color{purple}{M^{-1}}\color{blue}{A^T_y}+
  \color{blue}{A_t}\color{purple}{M}^{-1}\color{blue}{A_t^T}+
  \color{blue}{A_x}\color{purple}{M}^{-1}\color{blue}{A_x^T}) & (\color{blue}{A_t}\color{purple}{M^{-1}M^T_{\phi\sigma}}+\color{blue}{A_x}\color{purple}{M}^{-1}\color{cadmiumgreen}{L_x}
  +\color{blue}{A_y}\color{purple}{M}^{-1}\color{cadmiumgreen}{L_y}
  ) \\
  0 & 0 & -\color{cadmiumgreen}{L_y^T} & \color{cadmiumgreen}{L_x^T} \color{purple}{M^{-1}}\color{blue}{A_x^T}
 +\color{purple}{M_{\phi\sigma}}\color{purple}{M}^{-1}\color{blue}{A_t^T} & 
 -\color{cadmiumgreen}{L_x^T}\color{purple}{M^{-1}}\color{cadmiumgreen}{L_x}
  -\color{purple}{M_{\phi\sigma} M^{-1}M^T_{\phi\sigma}}+\color{purple}{M_{\phi}} + \color{red}{K} \\  \end{bmatrix} \notag \\
    \begin{bmatrix}
  \sigma_* \\ \sigma_x \\ \sigma_y \\ \lambda \\ \phi 
  \end{bmatrix}
  =
  \begin{bmatrix}
  0 \\ 0 \\ 0 \\ f \\ 0
  \end{bmatrix} \notag \\
  \end{eqnarray}      }}
      \item We add the third row multiplied by $\color{cadmiumgreen}{L^T_y}$ to the fifth row

{\footnotesize{\begin{eqnarray}
  \begin{bmatrix}
  {I} & 0 & 0 & \color{purple}{M}^{-1}\color{blue}{A_t^T} & -\color{purple}{M^{-1}M^T_{\phi\sigma}} \\    0 & I & 0 & \color{purple}{M^{-1}}\color{blue}{A_x^T} & -\color{purple}{M^{-1}}\color{cadmiumgreen}{L_x}  \\
  0 & 0 & I & \color{purple}{M^{-1}}\color{blue}{A_y^T} & -\color{purple}{M^{-1}}\color{cadmiumgreen}{L_y} \\
  0 & 0 & 0 & 
  -(\color{blue}{A_y}\color{purple}{M^{-1}}\color{blue}{A^T_y}+
  \color{blue}{A_t}\color{purple}{M}^{-1}\color{blue}{A_t^T}+
  \color{blue}{A_x}\color{purple}{M}^{-1}\color{blue}{A_x^T}) & (\color{blue}{A_t}\color{purple}{M^{-1}M^T_{\phi\sigma}}+\color{blue}{A_x}\color{purple}{M}^{-1}\color{cadmiumgreen}{L_x}
  +\color{blue}{A_y}\color{purple}{M}^{-1}\color{cadmiumgreen}{L_y}
  ) \\
 0 & 0 & 0 & \color{cadmiumgreen}{L_y^T} \color{purple}{M^{-1}}\color{blue}{A_y^T}
 + \color{cadmiumgreen}{L_x^T} \color{purple}{M^{-1}}\color{blue}{A_x^T}
  +\color{purple}{M_{\phi\sigma}}\color{purple}{M}^{-1}\color{blue}{A_t^T} & 
 \left(-\color{cadmiumgreen}{L_y^T}\color{purple}{M^{-1}}\color{cadmiumgreen}{L_y}
 -\color{cadmiumgreen}{L_x^T}\color{purple}{M^{-1}}\color{cadmiumgreen}{L_x}
   + \color{red}{K}\right)+  \\  
   &  &  &  & 
  \left(-\color{purple}{M_{\phi\sigma} M^{-1}M^T_{\phi\sigma}}+\color{purple}{M_{\phi}}\right)  \\   \end{bmatrix} \notag \\
    \begin{bmatrix}
  \sigma_* \\ \sigma_x \\ \sigma_y \\ \lambda \\ \phi 
  \end{bmatrix}
  =
  \begin{bmatrix}
  0 \\ 0 \\ 0 \\ f \\ 0
  \end{bmatrix} \notag  \\
  \end{eqnarray} }}     
  \item We notice that $-\color{cadmiumgreen}{L_y^T}\color{purple}{M^{-1}}\color{cadmiumgreen}{L_y}
  -\color{cadmiumgreen}{L_x^T}\color{purple}{M^{-1}}\color{cadmiumgreen}{L_x}+\color{red}{K}=\color{black}{0}$. 

This can be demonstrated as follows: 
Given a function~$u_h$ in some discrete space, let~$[u_h]$ denote the vector of its coefficients
in a fixed basis of this discrete space.
Then, using the connection between the relevant matrices
and their associated bilinear forms, we have for all~$u_h$, $w_h$
\begin{equation}
  [w_h]^T M [u_h] = \Prod{u_h}{w_h}, \qquad 
  [w_h]^T L_\alpha [u_h] = \Prod{\mathcal{L}_\alpha u_h}{w_h},
\end{equation}
for~$\alpha = x, y$.
Combining these, we see that
\begin{equation}
  [w_h]^T M [\Lop_\alpha u_h] = \Prod{\Lop_\alpha u_h}{w_h} = [w_h] L_\alpha [u_h]
  \Rightarrow
  M^{-1} L_\alpha [u_h] = [\Lop_\alpha u_h]
\end{equation}
and thus
\begin{equation}
  [w_h]^T L_\alpha^T M^{-1} L_\alpha [u_h] =
  [w_h]^T L_\alpha^T [\Lop_\alpha u_h] =
  \Prod{\Lop_\alpha u_h}{\Lop_\alpha w_h}.
\end{equation}
Finally,
\begin{align*}
  [w_h]^T (L_x^T M^{-1} L_x +  L_y^T M^{-1} L_y) [u_h] &=
  \Prod{\Lop_x u_h}{\Lop_x w_h} + \Prod{\Lop_y u_h}{\Lop_y w_h}\\ &=
  [w_h]^T K [u_h],
\end{align*}
by definition of~$K$.
Since~$u_h$, $w_h$ are arbitrary, $K = L_x^T M^{-1} L_x +  L_y^T M^{-1} L_y$.

{\footnotesize{\begin{eqnarray}
  \begin{bmatrix}
  {I} & 0 & 0 & \color{purple}{M}^{-1}\color{blue}{A_t^T} & -\color{purple}{M^{-1}M^T_{\phi\sigma}} \\    0 & I & 0 & \color{purple}{M^{-1}}\color{blue}{A_x^T} & -\color{purple}{M^{-1}}\color{cadmiumgreen}{L_x}  \\
  0 & 0 & I & \color{purple}{M^{-1}}\color{blue}{A_y^T} & -\color{purple}{M^{-1}}\color{cadmiumgreen}{L_y} \\
  0 & 0 & 0 & 
  -(\color{blue}{A_y}\color{purple}{M^{-1}}\color{blue}{A^T_y}+
  \color{blue}{A_t}\color{purple}{M}^{-1}\color{blue}{A_t^T}+
  \color{blue}{A_x}\color{purple}{M}^{-1}\color{blue}{A_x^T}) & (\color{blue}{A_t}\color{purple}{M^{-1}M^T_{\phi\sigma}}+\color{blue}{A_x}\color{purple}{M}^{-1}\color{cadmiumgreen}{L_x}
  +\color{blue}{A_y}\color{purple}{M}^{-1}\color{cadmiumgreen}{L_y}
  ) \\
 0 & 0 & 0 & \color{cadmiumgreen}{L_y^T} \color{purple}{M^{-1}}\color{blue}{A_y^T}
 + \color{cadmiumgreen}{L_x^T} \color{purple}{M^{-1}}\color{blue}{A_x^T}
  +\color{purple}{M_{\phi\sigma}}\color{purple}{M}^{-1}\color{blue}{A_t^T} & 
  -\color{purple}{M_{\phi\sigma} M^{-1}M^T_{\phi\sigma}}+\color{purple}{M_{\phi}}  \\   \end{bmatrix} \notag \\
    \begin{bmatrix}
  \sigma_* \\ \sigma_x \\ \sigma_y \\ \lambda \\ \phi 
  \end{bmatrix}
  =
  \begin{bmatrix}
  0 \\ 0 \\ 0 \\ f \\ 0
  \end{bmatrix} \notag  \\
  \end{eqnarray}}}      

  \item We notice that $(
  \color{blue}{A_x}\color{purple}{M}^{-1}\color{blue}{A_x^T}+\color{blue}{A_y}\color{purple}{M^{-1}}\color{blue}{A^T_y}+
  \color{blue}{A_t}\color{purple}{M}^{-1}\color{blue}{A_t^T})\color{black}{=S}$ is equal to the space-time stiffness matrix $S$. 

The argument is similar to the one in the previous point.
By definition of~$A_\alpha$ and using integration by parts we get
\begin{align*}
  [w_h]^T A_\alpha [u_h] = \Prod{D_\alpha u_h}{w_h} = -\Prod{u_h}{D_\alpha w_h} = - [w_h]^T A_\alpha^T [u_h].
\end{align*}
Therefore
\begin{align*}
  [w_h]^T M [D_\alpha u_h] = \Prod{D_\alpha u_h}{w_h} = [w_h] A_\alpha [u_h]
\end{align*}
and
\begin{equation*}
  [w_h]^T M [D_\alpha u_h] = \Prod{D_\alpha u_h}{w_h} = -\Prod{u_h}{D_\alpha w_h} = - [w_h] A_\alpha^T [u_h].
\end{equation*}
The above shows that
\begin{equation}
  [D_\alpha u_h] = M^{-1} A_\alpha [u_h] = -M^{-1} A_\alpha^T [u_h]
\end{equation}
and so
\begin{equation}
  [w_h]^T A_\alpha M^{-1} A_\alpha^T [u_h] = - [w_h]^T A_\alpha [D_\alpha u_h]
  = [w_h]^T A_\alpha^T [D_\alpha u_h] = \Prod{D_\alpha w_h}{D_\alpha u_h}
\end{equation}
Summing up these equations for each direction we get
\begin{equation}
[w_h]^T (A_x M^{-1} A_x^T + A_y M^{-1} A_y^T + A_t M^{-1} A_t^T) [u_h] =
\Prod{\nabla_{x,t} u_h}{\nabla_{x,t} w_h} = [w_h]^T S [u_h]
\end{equation}
  
{\footnotesize{\begin{eqnarray}
  \begin{bmatrix}
  {I} & 0 & 0 & \color{purple}{M}^{-1}\color{blue}{A_t^T} & -\color{purple}{M^{-1}M^T_{\phi\sigma}} \\    0 & I & 0 & \color{purple}{M^{-1}}\color{blue}{A_x^T} & -\color{purple}{M^{-1}}\color{cadmiumgreen}{L_x}  \\
  0 & 0 & I & \color{purple}{M^{-1}}\color{blue}{A_y^T} & -\color{purple}{M^{-1}}\color{cadmiumgreen}{L_y} \\
  0 & 0 & 0 & 
  -S & (\color{blue}{A_t}\color{purple}{M^{-1}M^T_{\phi\sigma}}+\color{blue}{A_x}\color{purple}{M}^{-1}\color{cadmiumgreen}{L_x}
  +\color{blue}{A_y}\color{purple}{M}^{-1}\color{cadmiumgreen}{L_y}
  ) \\
 0 & 0 & 0 & \color{cadmiumgreen}{L_y^T} \color{purple}{M^{-1}}\color{blue}{A_y^T}
 + \color{cadmiumgreen}{L_x^T} \color{purple}{M^{-1}}\color{blue}{A_x^T}
  +\color{purple}{M_{\phi\sigma}}\color{purple}{M}^{-1}\color{blue}{A_t^T} & 
  -\color{purple}{M_{\phi\sigma} M^{-1}M^T_{\phi\sigma}}+\color{purple}{M_{\phi}}  \\   \end{bmatrix} \notag \\
    \begin{bmatrix}
  \sigma_* \\ \sigma_x \\ \sigma_y \\ \lambda \\ \phi 
  \end{bmatrix}
  =
  \begin{bmatrix}
  0 \\ 0 \\ 0 \\ f \\ 0
  \end{bmatrix} \notag  \\
  \end{eqnarray}}}   
  We notice that~$M_{\phi\sigma}M^{-1}M_{\phi\sigma}^T = M_\phi$.
  This can be demonstrated as follows.
  Let~$V_h$ and~$U_h$ denote the discrete spaces of~$\phi$ and components of~$\Sig$.
  Since they differ only by the Dirichlet boundary conditions, we can order the basis functions
  so that~$V_h = \operatorname{span}\{e_1,\ldots,e_n\}$
  and~$U_h = \operatorname{span}\{e_1,\ldots,e_n,e_{n+1},\ldots,e_{N}\}$.
  Consequently, we can partition the above matrices (possibly after permuting rows and columns)
  into blocks corresponding to interior and boundary basis functions as
  \begin{equation}
    M =
    \begin{bmatrix}
      M_\phi & M_{\phi\partial} \\
      M_{\phi\partial}^T & M_\partial
    \end{bmatrix}
    \qquad
    M_{\phi\sigma} =
    \begin{bmatrix}
      M_\phi & M_{\phi\partial}
    \end{bmatrix}
  \end{equation}
  where~$M_{\phi\partial}$ and~$M_\partial$ contain~$L^2$ products of boundary basis
  functions with interior basis functions and boundary basis functions, respectively.
  More concretely, $(M_\phi)_{ij} = \Prod{e_j}{e_i}$ for~$i, j = 1,\ldots,n$,
  $(M_\partial)_{ij} = \Prod{e_{n+j}}{e_{n+i}}$ for~$i, j = 1,\ldots,N - n$
  and~$(M_{\phi\partial})_{ij} = \Prod{e_{n+j}}{e_i}$ for~$i = 1,\ldots,n$,
  $j = 1,\ldots,N-n$.
  Writing~$M^{-1}$ as
  \begin{equation}
    M^{-1} =
    \begin{bmatrix}
      A & B \\ C & D  
    \end{bmatrix}
  \end{equation}
  with the same block layout as for~$M$,
  we have by definition of the inverse~$A M_\phi + B M_{\phi\partial}^T = I$
  and~$C M_\phi + D M_{\phi\partial}^T = 0$, thus
  \begin{equation}
    M^{-1} M_{\phi\sigma}^T =
    \begin{bmatrix}
      A & B \\ C & D  
    \end{bmatrix}
    \begin{bmatrix}
      M_\phi \\ M_{\phi\partial}^T
    \end{bmatrix}
    =
    \begin{bmatrix}
      I \\ 0
    \end{bmatrix}
  \end{equation}
  and finally
  \begin{equation}
    M_{\phi\sigma} M^{-1} M_{\phi\sigma}^T =
    \begin{bmatrix}
      M_\phi & M_{\phi\partial}
    \end{bmatrix}
    \begin{bmatrix}
      I \\ 0
    \end{bmatrix}
    =
    M_\phi
  \end{equation}
  We end up with
{\footnotesize{  \begin{eqnarray}
  \begin{bmatrix}
  {I} & 0 & 0 & \color{purple}{M}^{-1}\color{blue}{A_t^T} & -\color{purple}{M^{-1}M^T_{\phi\sigma}} \\    0 & I & 0 & \color{purple}{M^{-1}}\color{blue}{A_x^T} & -\color{purple}{M^{-1}}\color{cadmiumgreen}{L_x}  \\
  0 & 0 & I & \color{purple}{M^{-1}}\color{blue}{A_y^T} & -\color{purple}{M^{-1}}\color{cadmiumgreen}{L_y} \\
  0 & 0 & 0 & 
  -S & (\color{blue}{A_t}\color{purple}{M^{-1}M^T_{\phi\sigma}}+\color{blue}{A_x}\color{purple}{M}^{-1}\color{cadmiumgreen}{L_x}
  +\color{blue}{A_y}\color{purple}{M}^{-1}\color{cadmiumgreen}{L_y}
  ) \\
 0 & 0 & 0 & \color{cadmiumgreen}{L_y^T} \color{purple}{M^{-1}}\color{blue}{A_y^T}
 + \color{cadmiumgreen}{L_x^T} \color{purple}{M^{-1}}\color{blue}{A_x^T}
  +\color{purple}{M_{\phi\sigma}}\color{purple}{M}^{-1}\color{blue}{A_t^T} & 
  0  \\   \end{bmatrix} \notag \\
    \begin{bmatrix}
  \sigma_* \\ \sigma_x \\ \sigma_y \\ \lambda \\ \phi 
  \end{bmatrix}
  =
  \begin{bmatrix}
  0 \\ 0 \\ 0 \\ f \\ 0
  \end{bmatrix} \notag  \\
  \end{eqnarray}}}   

\item  In the reminding system 

  {\footnotesize{\begin{eqnarray}
  \begin{bmatrix}
  -S & (\color{blue}{A_t}\color{purple}{M^{-1}M^T_{\phi\sigma}}+\color{blue}{A_x}\color{purple}{M}^{-1}\color{cadmiumgreen}{L_x}
  +\color{blue}{A_y}\color{purple}{M}^{-1}\color{cadmiumgreen}{L_y}
  ) \\
  (\color{cadmiumgreen}{L_y^T} \color{purple}{M^{-1}}\color{blue}{A_y^T}
 + \color{cadmiumgreen}{L_x^T} \color{purple}{M^{-1}}\color{blue}{A_x^T}
  +\color{purple}{M_{\phi\sigma}}\color{purple}{M}^{-1}\color{blue}{A_t^T}) & 
  0  \\   \end{bmatrix} 
    \begin{bmatrix}
   \lambda \\ \phi 
  \end{bmatrix}
  =
  \begin{bmatrix}
  f \\ 0
  \end{bmatrix} \notag \\
  \end{eqnarray} }}  

we denote
  $Z=(\color{blue}{A_t}\color{purple}{M^{-1}M^T_{\phi\sigma}}+\color{blue}{A_x}\color{purple}{M}^{-1}\color{cadmiumgreen}{L_x}
  +\color{blue}{A_y}\color{purple}{M}^{-1}\color{cadmiumgreen}{L_y}
  )$ to obtain

  \begin{eqnarray}
  \begin{bmatrix}
  -S & Z \\
  Z^T & 
  0  \\   \end{bmatrix} 
    \begin{bmatrix}
   \lambda \\ \phi 
  \end{bmatrix}
  =
  \begin{bmatrix}
  f \\ 0
  \end{bmatrix}
  \end{eqnarray}   

  and $ZS^{-1}Z^T$ is symmetric and positive definite.
This system is of the Uzawa type \cite{Uzawa} and it can be solved efficiently using an iterative solver. We employ the GMRES solver summarized in Algorithm 1.

\begin{algorithm}
\caption{The GMRES algorithm} \label{Alg4}
\begin{algorithmic}[1]
\REQUIRE $A$ matrix, $b$ right-hand-side vector, $x_0$ starting point
\STATE Compute $r_0=b-Ax_0$
\STATE Compute $v_1=\frac{r_0}{\|r_0\|}$
\STATE {\bf for}
$j = 1, 2, . . . , k$
\STATE 
Compute $h_{i,j} = \left(Av_j , v_i\right)$
 for $i = 1, 2, . . . , j$
\STATE Compute $\hat{v}_{j+1}=Av_j-\sum_{i=1,...,j}h_{i,j}v_i$
\STATE Compute $h_{j+1,j} = \|\hat{v}_{j+1}\|_2$
\STATE Compute $v_{j+1}=\hat{v}_{j+1}/h_{j+1,j}$
\STATE {\bf end for}
\STATE Form solution $x_k=x_0+V_ky_k$, where 
$V_k=[v_1 ... v_k]$,
and $y_k$ minimizes
$J(y)=\|\beta e_1-\hat{H}_ky\|$, where
$\hat{H}=\begin{bmatrix}h_{1,1} & h_{1,2} \cdots h_{1,k}\\
h_{2,1} & h_{2,2} \cdots h_{2,k}\\
0 &  \ddots & \ddots & \vdots \\
\vdots & \ddots & h_{k,k-1} & h_{k,k} \\
0 & \cdots & 0 & h_{k+1,k}
\end{bmatrix}$.
\end{algorithmic}
\end{algorithm}
\end{enumerate}
  
\section{Numerical results with stabilization}\label{Stabilization}

\begin{figure}[h]
\centering
  \begin{subfigure}[b]{0.32\textwidth}
  \centering
  \includegraphics[width=\textwidth]{new-eps1e-3-s1.0-arrow}
  \end{subfigure}
  \begin{subfigure}[b]{0.32\textwidth}
  \centering
  \includegraphics[width=\textwidth]{new-eps1e-5-s1.0-arrow}
  \end{subfigure}
  \begin{subfigure}[b]{0.32\textwidth}
  \centering
  \includegraphics[width=\textwidth]{new-eps1e-6-s1.0-arrow}
  \end{subfigure}

  \begin{subfigure}[b]{0.32\textwidth}
  \includegraphics[width=\textwidth]{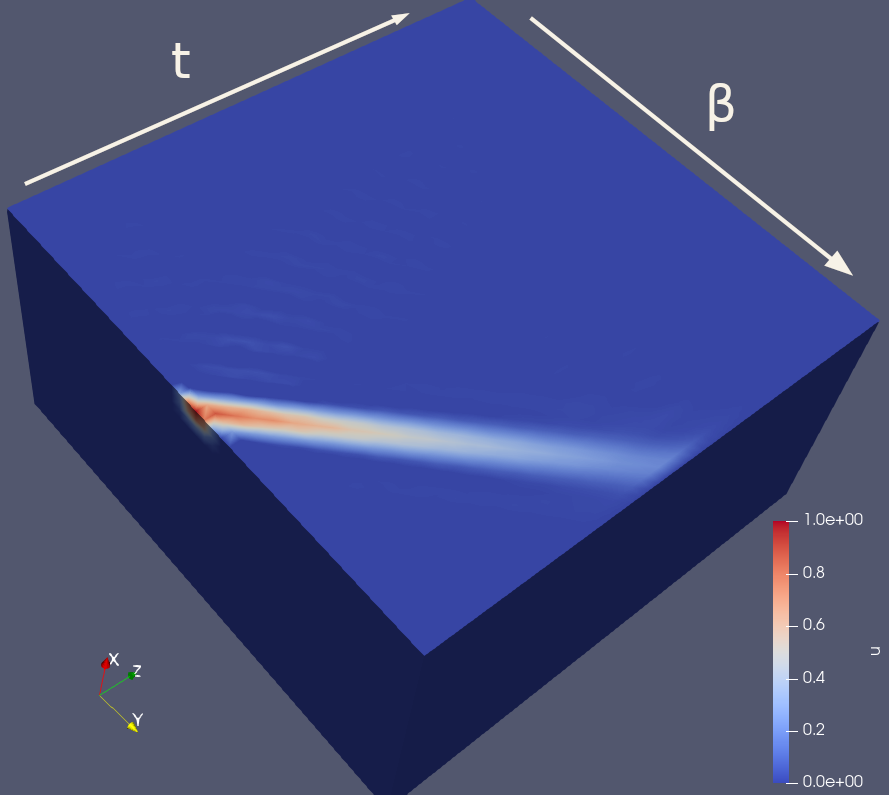}
  \caption{$\varepsilon = 10^{-3}$, $s = 1$}
  \end{subfigure}
  \begin{subfigure}[b]{0.32\textwidth}
  \centering
  \includegraphics[width=\textwidth]{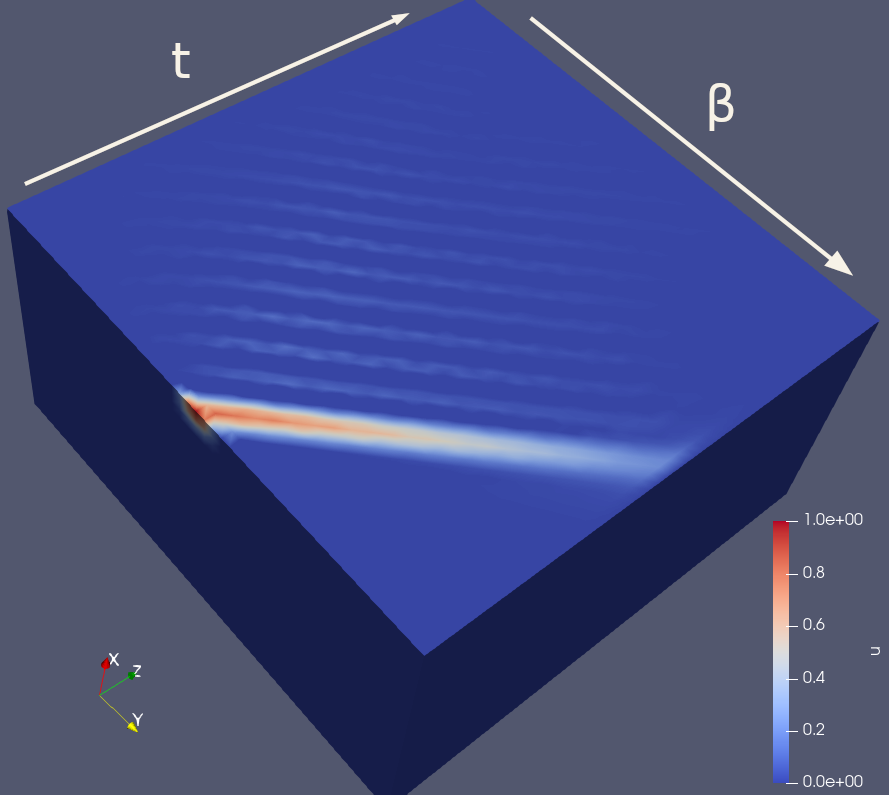}
  \caption{$\varepsilon = 10^{-5}$, $s = 0.5$}
  \end{subfigure}
  \begin{subfigure}[b]{0.32\textwidth}
  \centering
  \includegraphics[width=\textwidth]{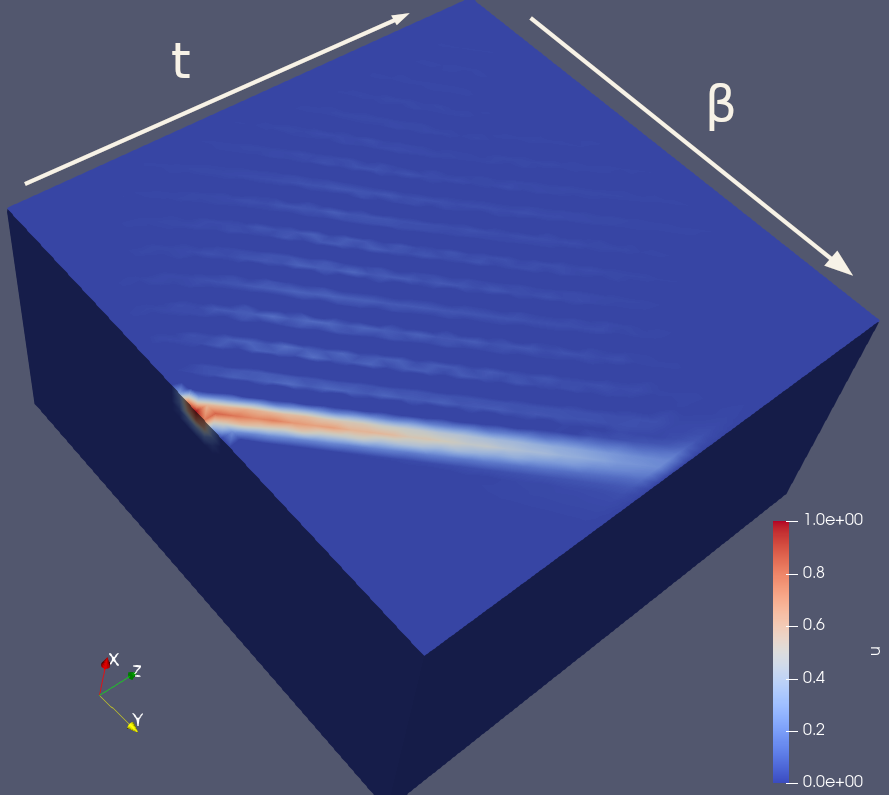}
  \caption{$\varepsilon = 10^{-6}$, $s = 1$}
  \end{subfigure}

  \caption{Results for~$\beta = (0, s)$, $T = 1$, a mesh with  $32\times 32 \times 32$ elements with quadratic B-splines (top: no stabilization, bottom: stabilization).}
  \label{fig:stab}
\end{figure}

We focus on our model problem
~\eqref{eq:adv-diff-problem}
formulated on a regular domain~$\Omega \times (0, T])= (0, 1)^3$ with $\boldsymbol{\beta} = (0, 0.3)^t$, $\varepsilon = 10^{-5}$, and 
$f = 0$, using the initial state $u_0$ from (\ref{eq:init}).
Now, we emply the stabilized formulation (\ref{eq:stab}).
The numerical results are summarized in Figure \ref{fig:stab}. We can read that our stabilization effort decreases the unwanted oscillations.
Additionally, in Table 1, we present the execution times for the solver with different dimensions of the computational grids.

\begin{table}[t]
\centering
\begin{tabular}{ccc}\hline 
  tensor product grid & DoFs &  solver [s] \\
  \hline
$16\times 16\times 16$ &    5,832 & 0.088 \\
$24\times 24\times 24$ &   13,824 & 0.679 \\
$32\times 32\times 32$ &   39,304 & 6.479 \\ 
$64\times 64\times 64$ &  287,496 & 54.758 \\ 
  \hline
\end{tabular}
\caption{Solver execution times for IGA grids}
\end{table}

\section{Numerical results with stabilization and adaptation}

\begin{figure}
  \centering
  \begin{subfigure}[b]{0.32\textwidth}
  \begin{tikzpicture}[scale=0.7]
  \node at (0,0){\includegraphics[width=1\textwidth]{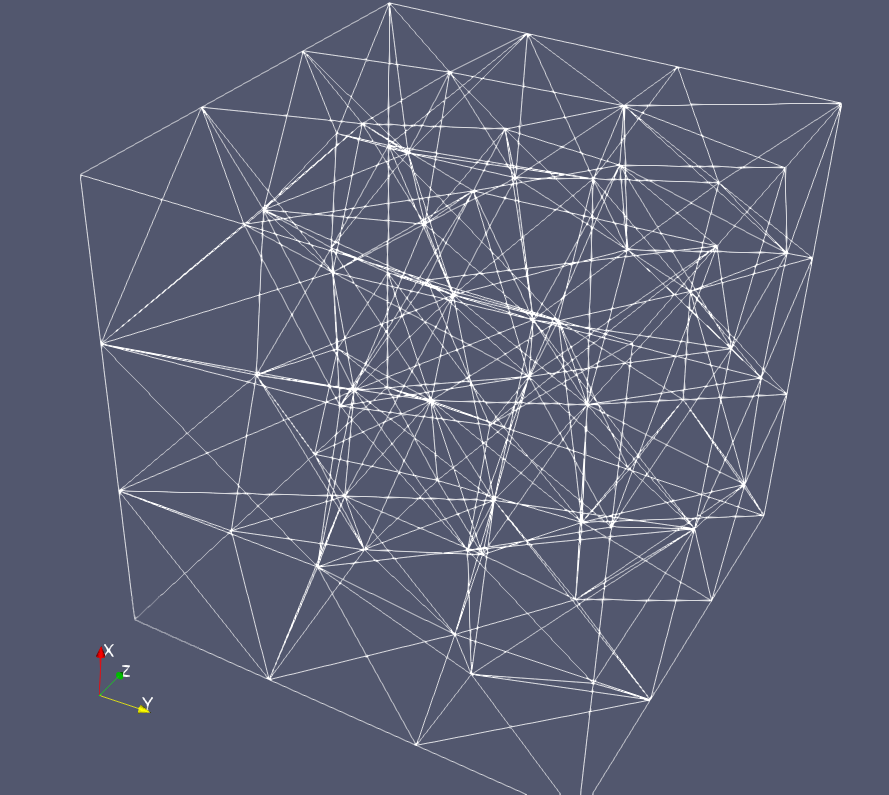}};
  \draw[white,->,thick] (-2.4,-1.6)--(-0.3,0.5);
  \node[white] at (-1.3,-1.){$t$};
  \end{tikzpicture} 
  \end{subfigure}
  \begin{subfigure}[b]{0.32\textwidth}
  \begin{tikzpicture}[scale=0.7]
  \node at (0,0)
  {\includegraphics[width=\textwidth]{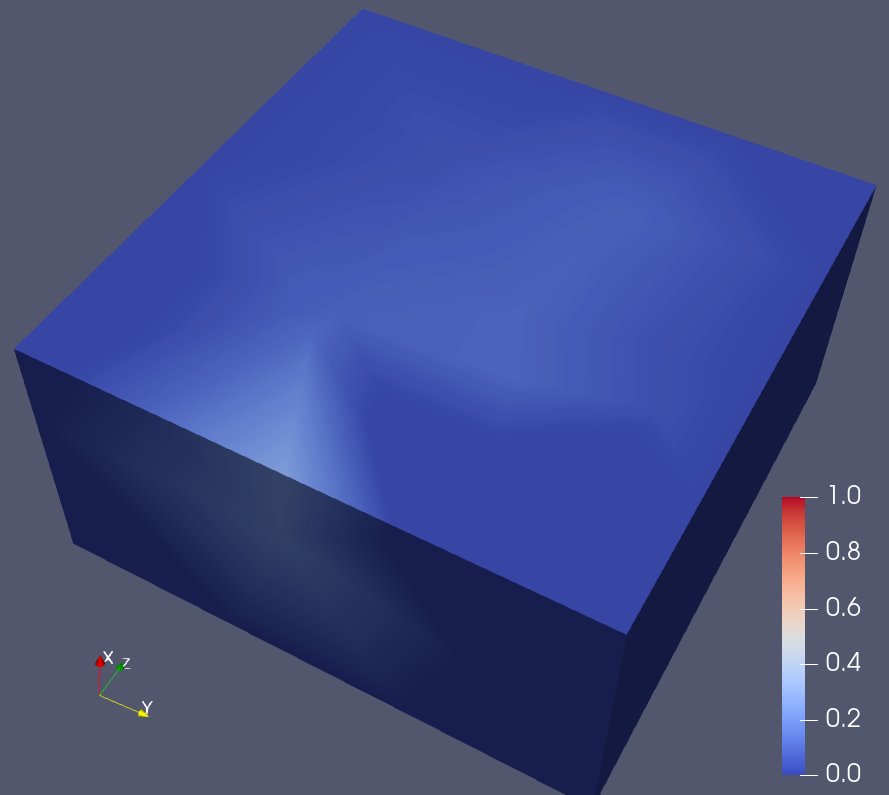} };
  \draw[white,->,thick] (-3.3,0.3)--(-1.2,2.4);
  \node[white] at (-2,2.2){$t$};
  \end{tikzpicture} 
  \end{subfigure}
  \begin{subfigure}[b]{0.32\textwidth}
  \begin{tikzpicture}[scale=0.7]
  \node at (0,0)
  {\includegraphics[width=1\textwidth]{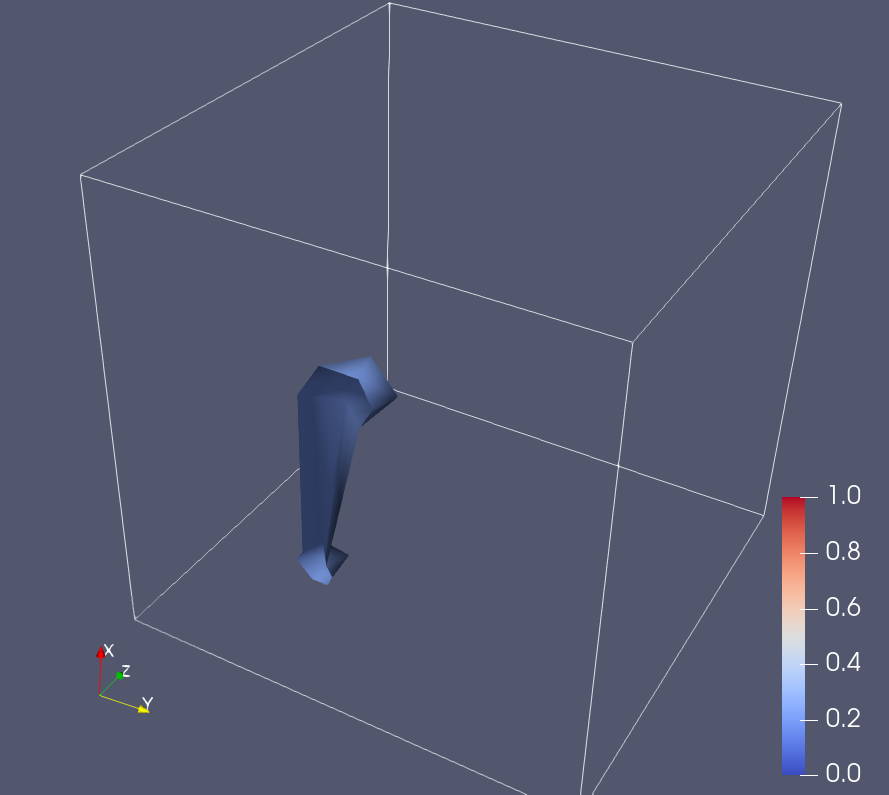} };
  \draw[white,->,thick] (-2.85,1.7)--(-1.2,2.7);
  \node[white] at (-2.2,2.6){$t$};
  \end{tikzpicture}
  \end{subfigure}
  \caption{Initial space-time mesh, the solution over the initial mesh, and the contour of the solution.}
  \label{fig:init-mesh}
\end{figure}
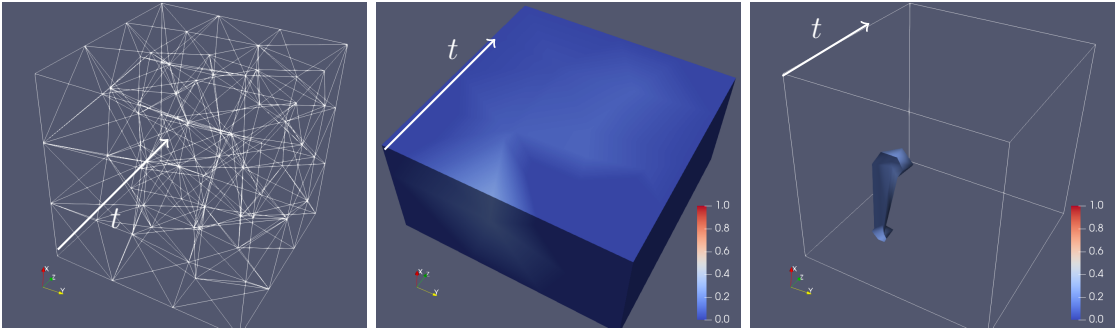

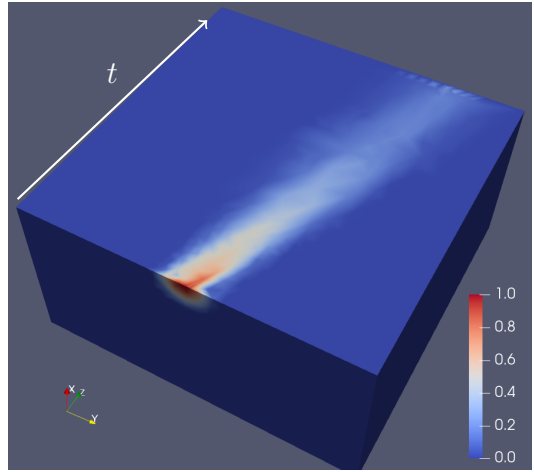
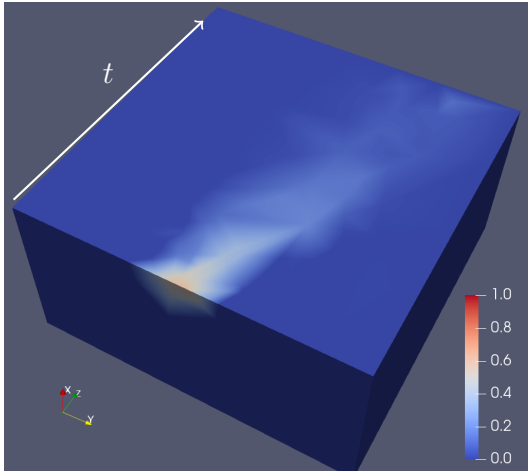
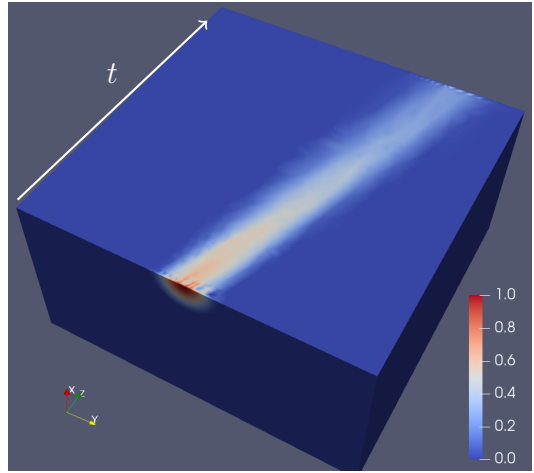
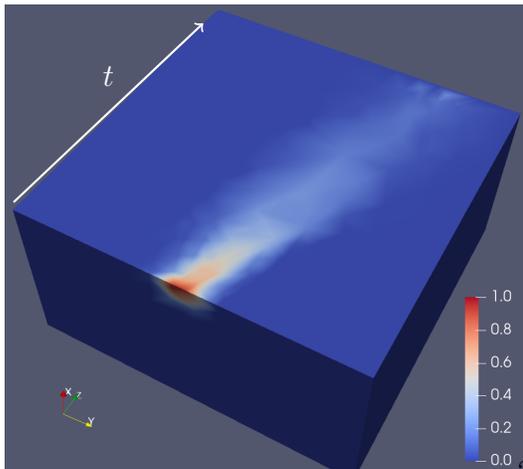
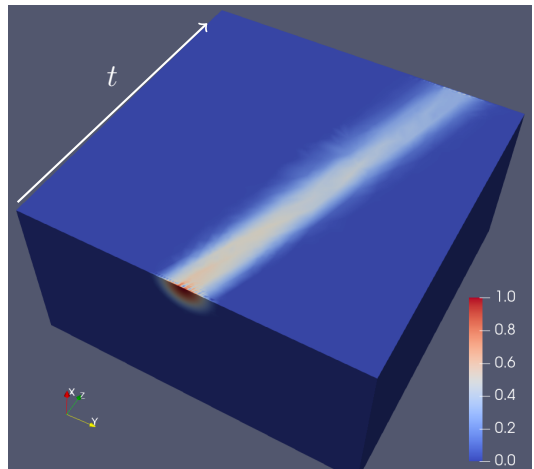
\begin{figure}
\centering
  \begin{subfigure}[b]{0.46\textwidth}
  \begin{tikzpicture}[scale=0.7]
  \node at (0,0)
  {\includegraphics[width=\textwidth]{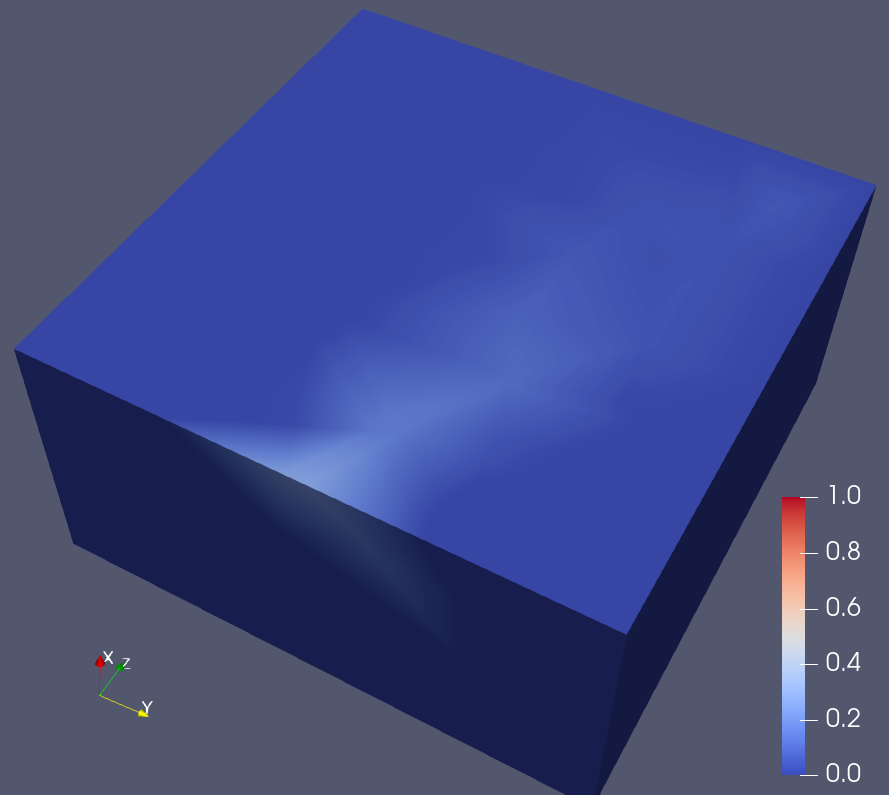} };
  \draw[white,->,thick] (-4.8,0.7)--(-1.2,4.1);
  \node[white] at (-3,3.1){$t$};
  \end{tikzpicture}
  \caption{refinement level 1}
  \end{subfigure}
  ~
  \begin{subfigure}[b]{0.46\textwidth}
  \begin{tikzpicture}[scale=0.7]
  \node at (0,0){
  \includegraphics[width=\textwidth]{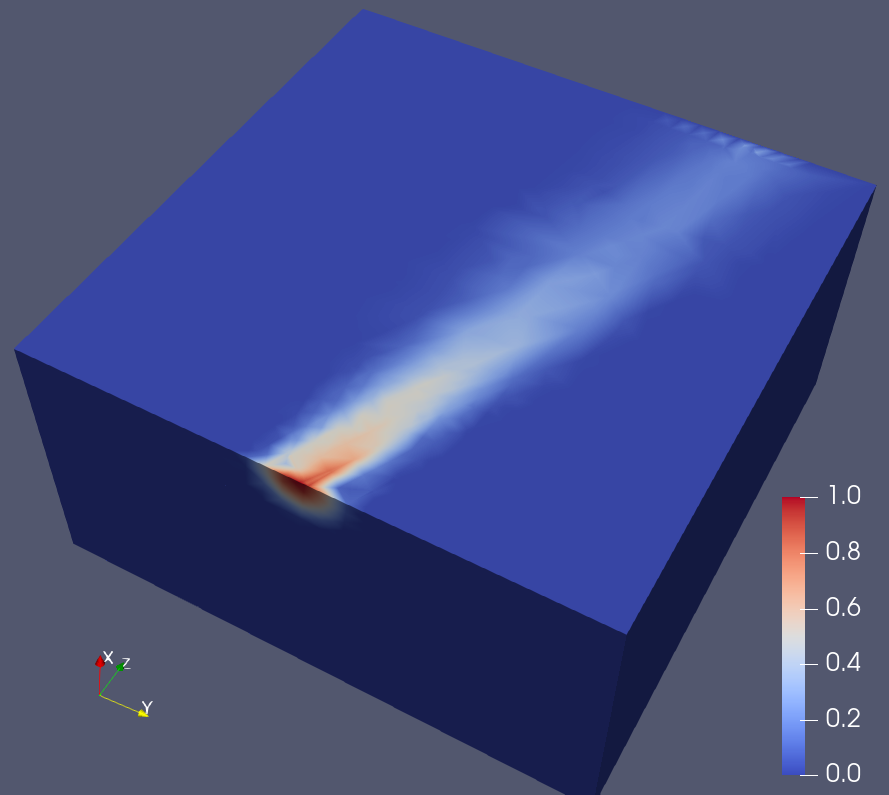}};
  \draw[white,->,thick] (-4.8,0.7)--(-1.2,4.1);
  \node[white] at (-3,3.1){$t$};
  \end{tikzpicture}
  \caption{refinement level 4}
  \end{subfigure}

  \begin{subfigure}[b]{0.46\textwidth}
  \begin{tikzpicture}[scale=0.7]
  \node at (0,0){
  \includegraphics[width=\textwidth]{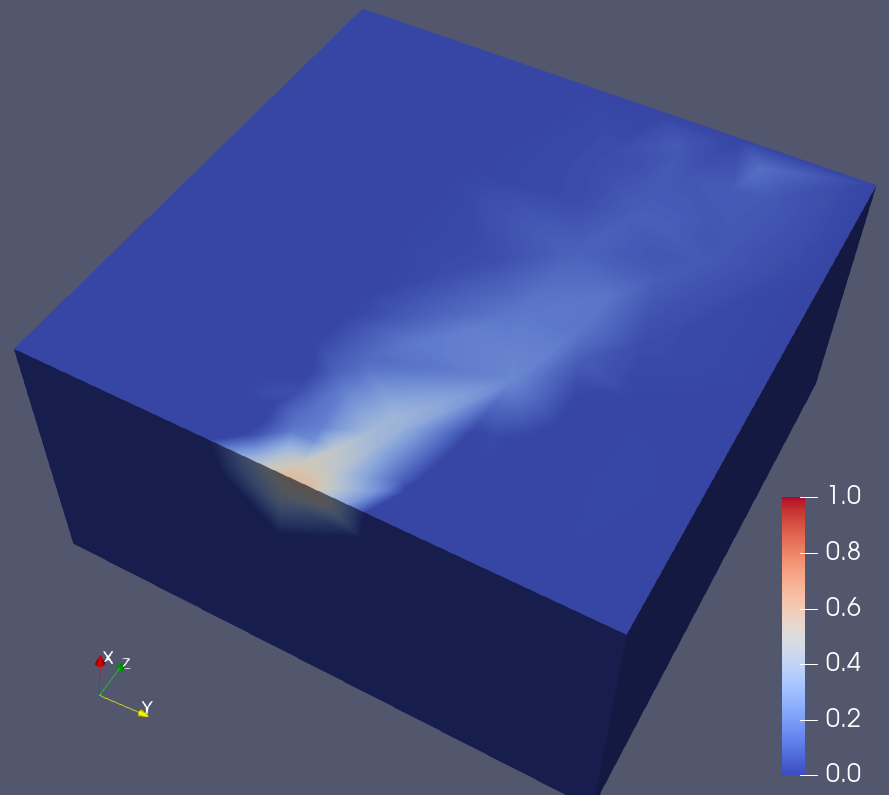}};
\draw[white,->,thick] (-4.8,0.7)--(-1.2,4.1);
  \node[white] at (-3,3.1){$t$};
  \end{tikzpicture}
  \caption{refinement level 2}
  \end{subfigure}
  ~
  \begin{subfigure}[b]{0.46\textwidth}
  \begin{tikzpicture}[scale=0.7]
  \node at (0,0){
  \includegraphics[width=\textwidth]{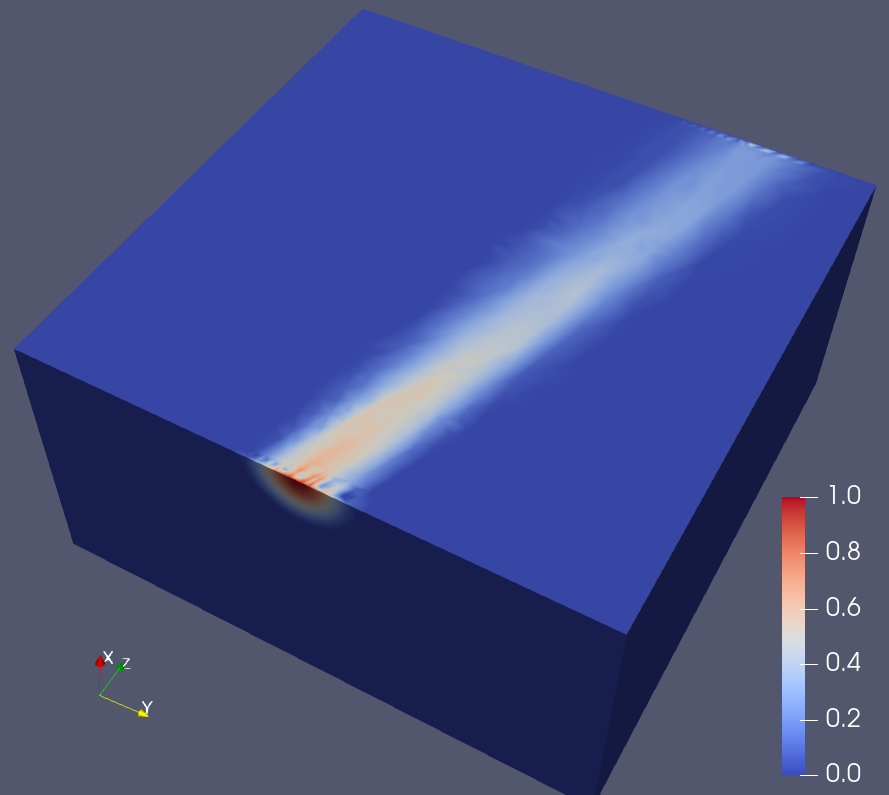}};
  \draw[white,->,thick] (-4.8,0.7)--(-1.2,4.1);
  \node[white] at (-3,3.1){$t$};
  \end{tikzpicture}
  \caption{refinement level 5}
  \end{subfigure}

  \begin{subfigure}[b]{0.46\textwidth}
  \begin{tikzpicture}[scale=0.7]
  \node at (0,0){
  \includegraphics[width=\textwidth]{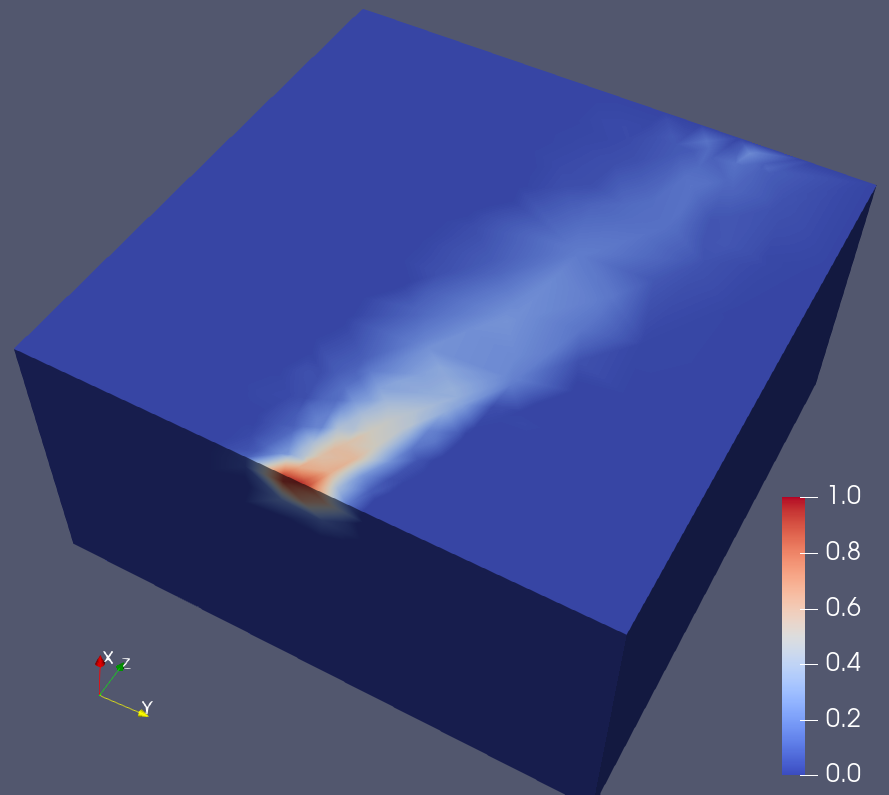}};
  \draw[white,->,thick] (-4.8,0.7)--(-1.2,4.1);
  \node[white] at (-3,3.1){$t$};
  \end{tikzpicture}
  \caption{refinement level 3}
  \end{subfigure}
  ~
  \begin{subfigure}[b]{0.46\textwidth}
  \begin{tikzpicture}[scale=0.7]
  \node at (0,0){
  \includegraphics[width=\textwidth]{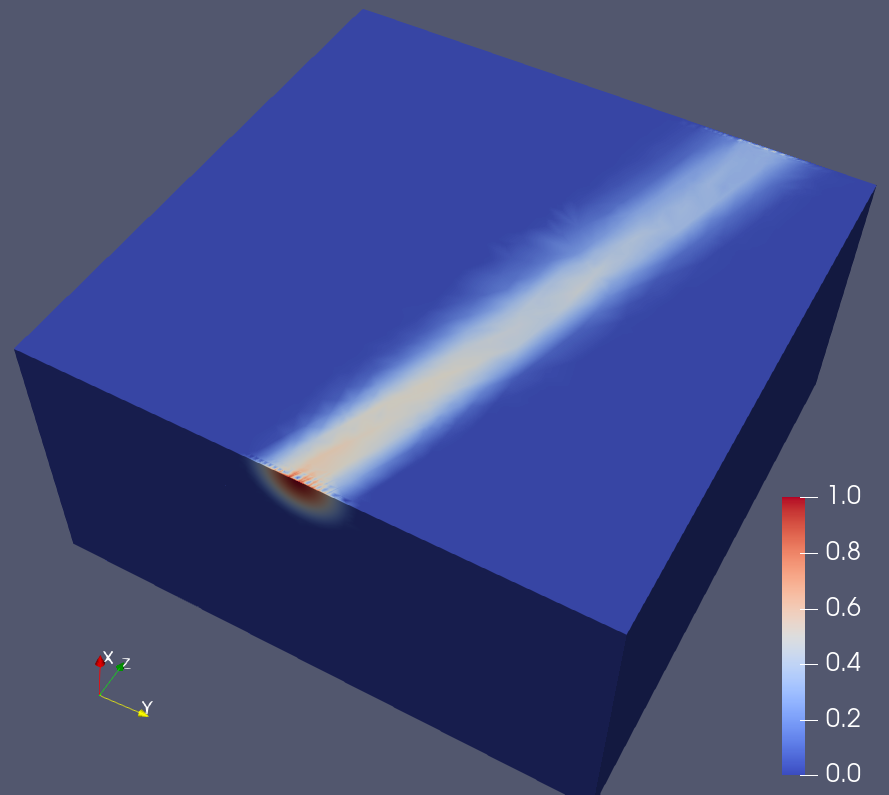} };
  \draw[white,->,thick] (-4.8,0.7)--(-1.2,4.1);
  \node[white] at (-3,3.1){$t$};
  \end{tikzpicture}
  \caption{refinement level 6}
  \end{subfigure}
  
  \caption{Solutions in subsequent space-time refinement steps}
\label{fig:adapt}
\end{figure}

\begin{figure}
\centering
  \begin{subfigure}[b]{0.46\textwidth}
  \centering
  \includegraphics[width=\textwidth]{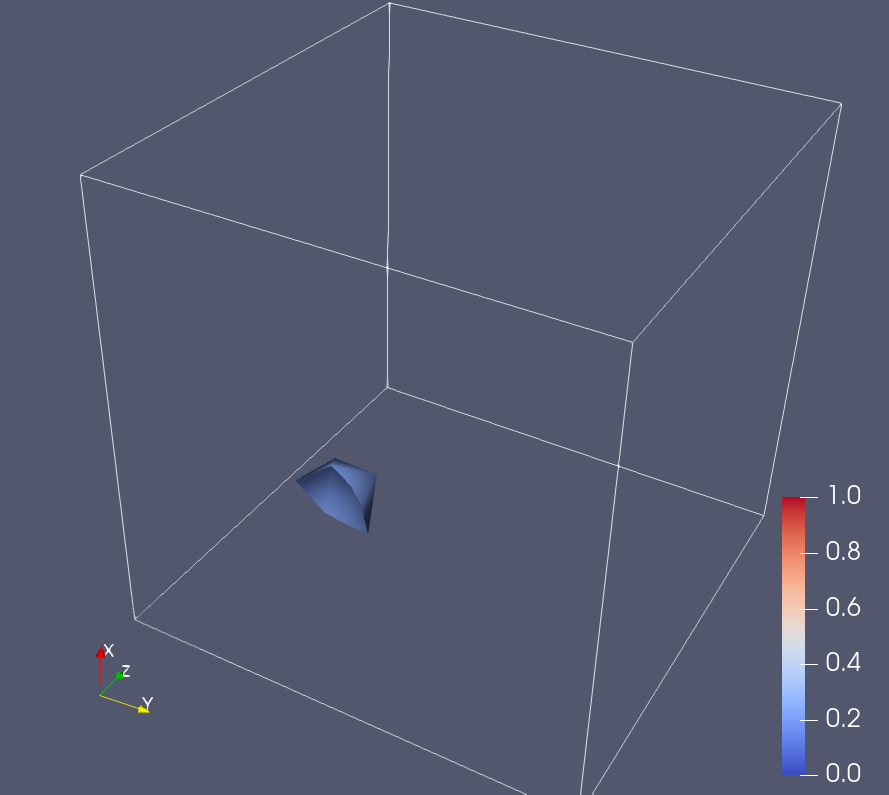} 
  \caption{refinement level 1}
  \end{subfigure}
  ~
  \begin{subfigure}[b]{0.46\textwidth}
  \centering
  \includegraphics[width=\textwidth]{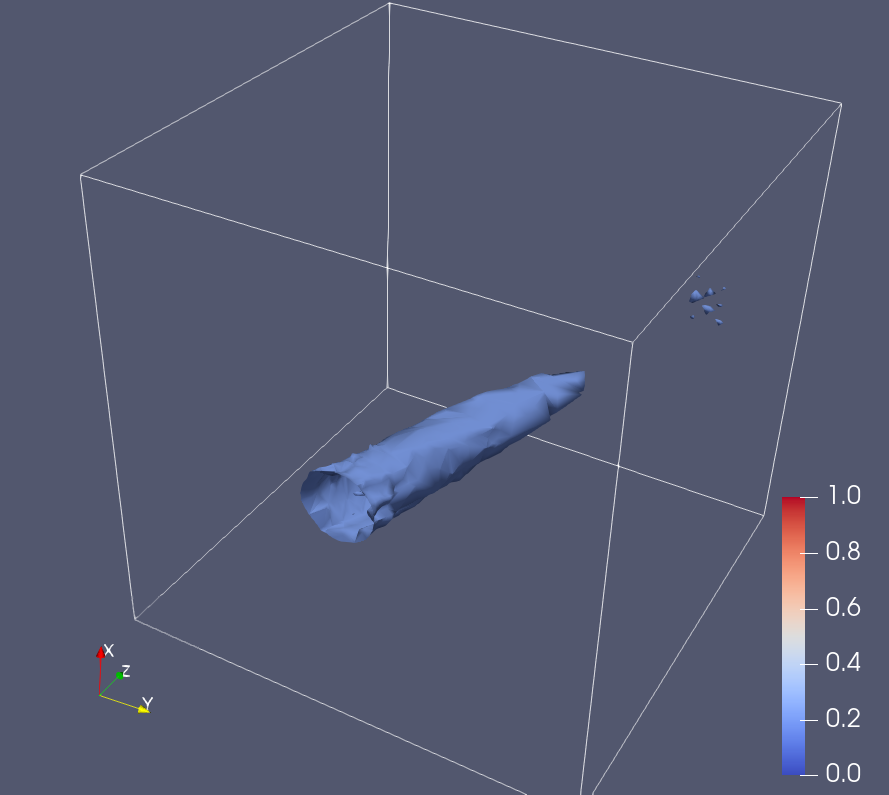} 
  \caption{refinement level 4}
  \end{subfigure}

  \begin{subfigure}[b]{0.46\textwidth}
  \centering
  \includegraphics[width=\textwidth]{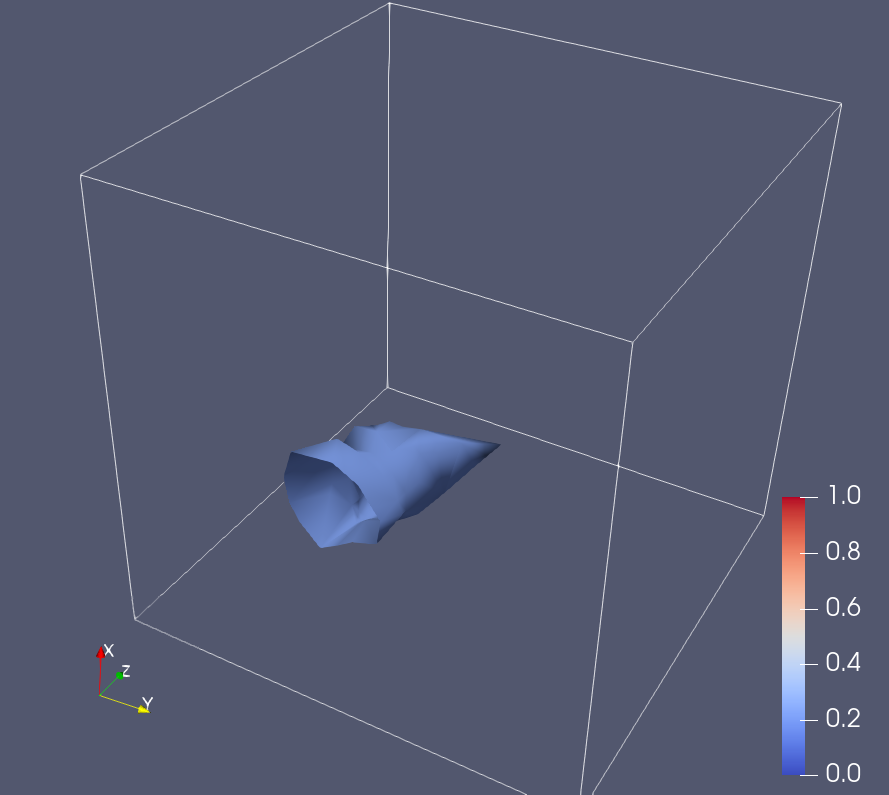} 
  \caption{refinement level 2}
  \end{subfigure}
  ~
  \begin{subfigure}[b]{0.46\textwidth}
  \centering
  \includegraphics[width=\textwidth]{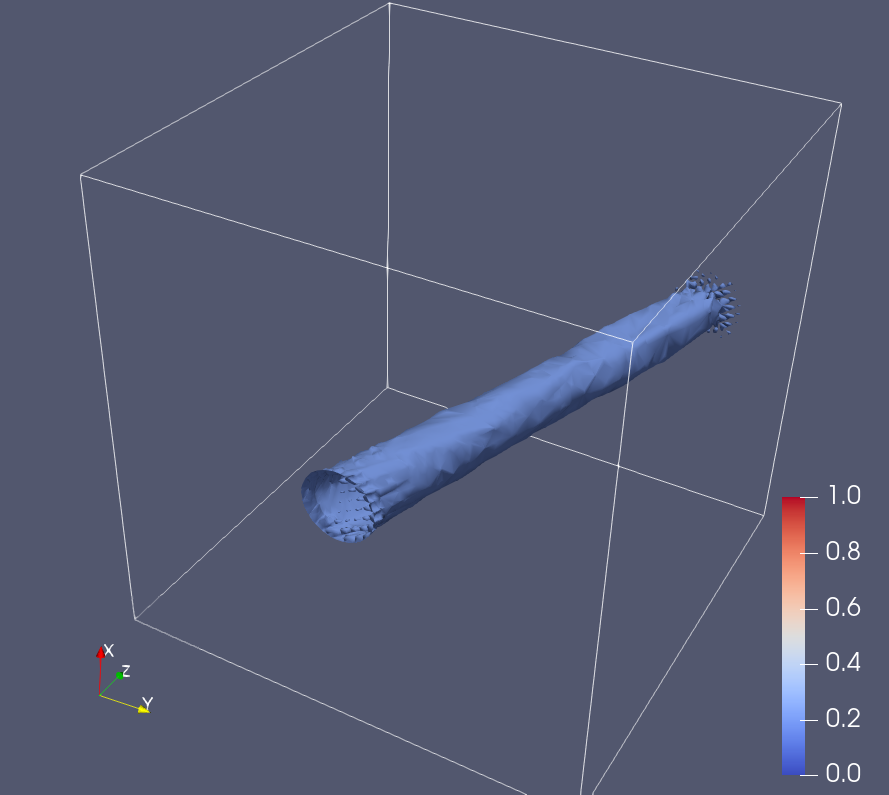} 
  \caption{refinement level 5}
  \end{subfigure}

  \begin{subfigure}[b]{0.46\textwidth}
  \centering
  \includegraphics[width=\textwidth]{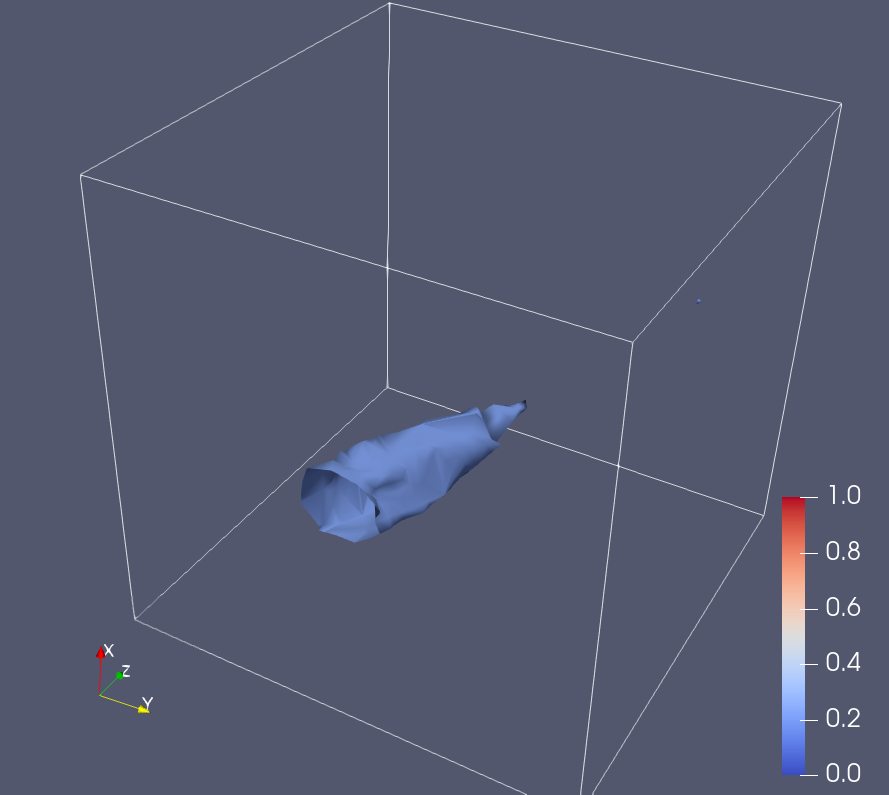} 
  \caption{refinement level 3}
  \end{subfigure}
  ~
  \begin{subfigure}[b]{0.46\textwidth}
  \centering
  \includegraphics[width=\textwidth]{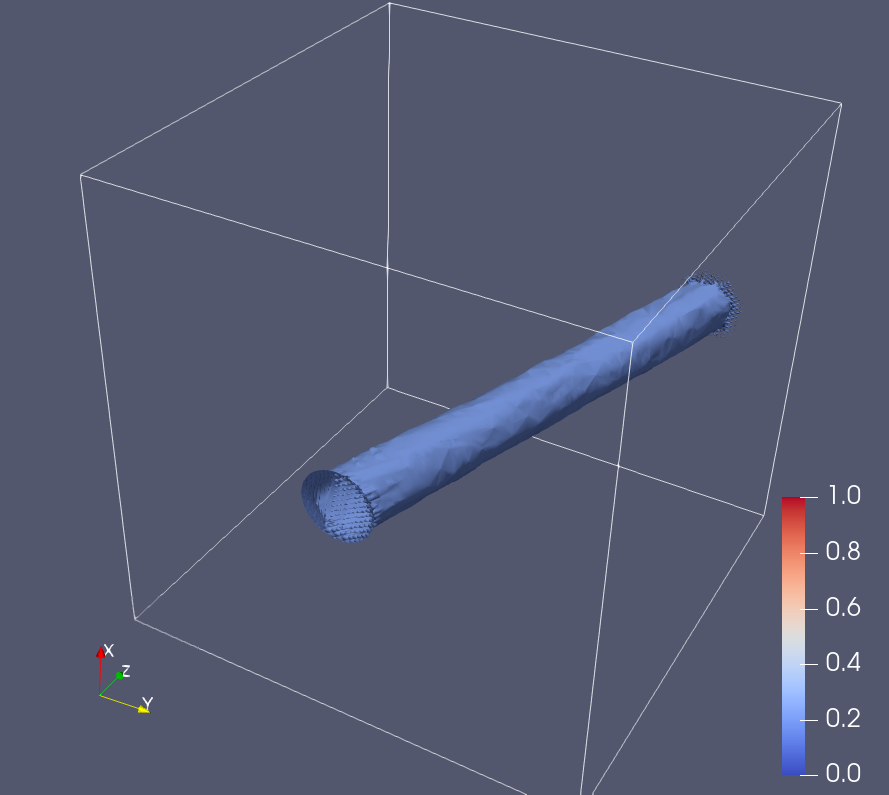} 
  \caption{refinement level 6}
  \end{subfigure}
  
  \caption{Contours ($\phi = 0.2$) in subsequent space-time refinement steps}
\label{fig:cont}
\end{figure}

\begin{figure}
\centering
  \begin{subfigure}[b]{0.46\textwidth}
  \centering
  \includegraphics[width=\textwidth]{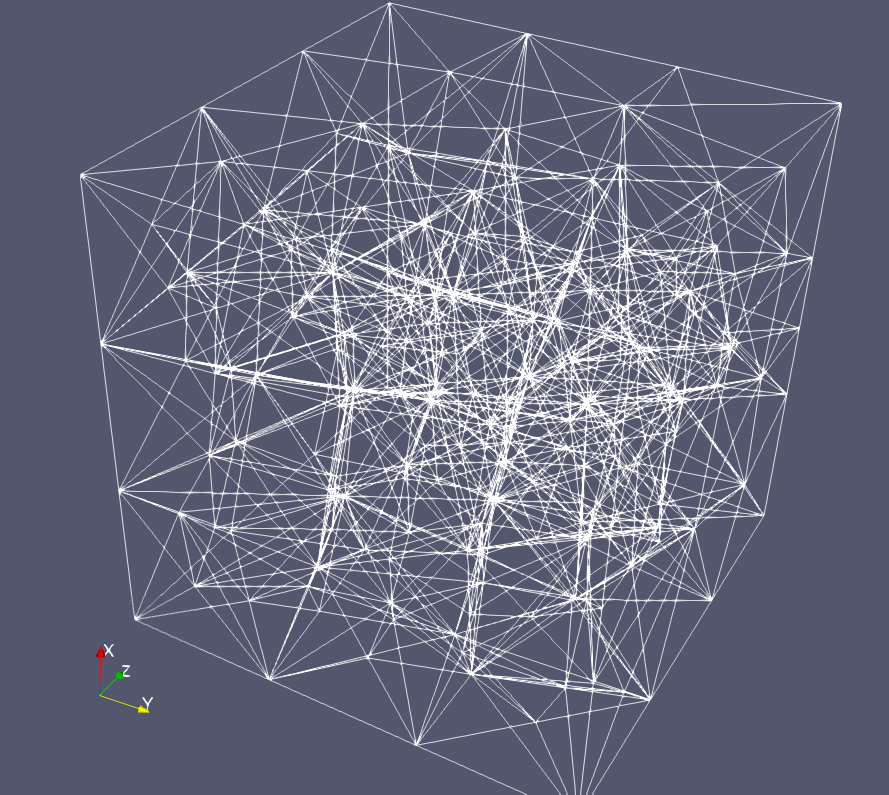} 
  \caption{refinement level 1}
  \end{subfigure}
  ~
  \begin{subfigure}[b]{0.46\textwidth}
  \centering
  \includegraphics[width=\textwidth]{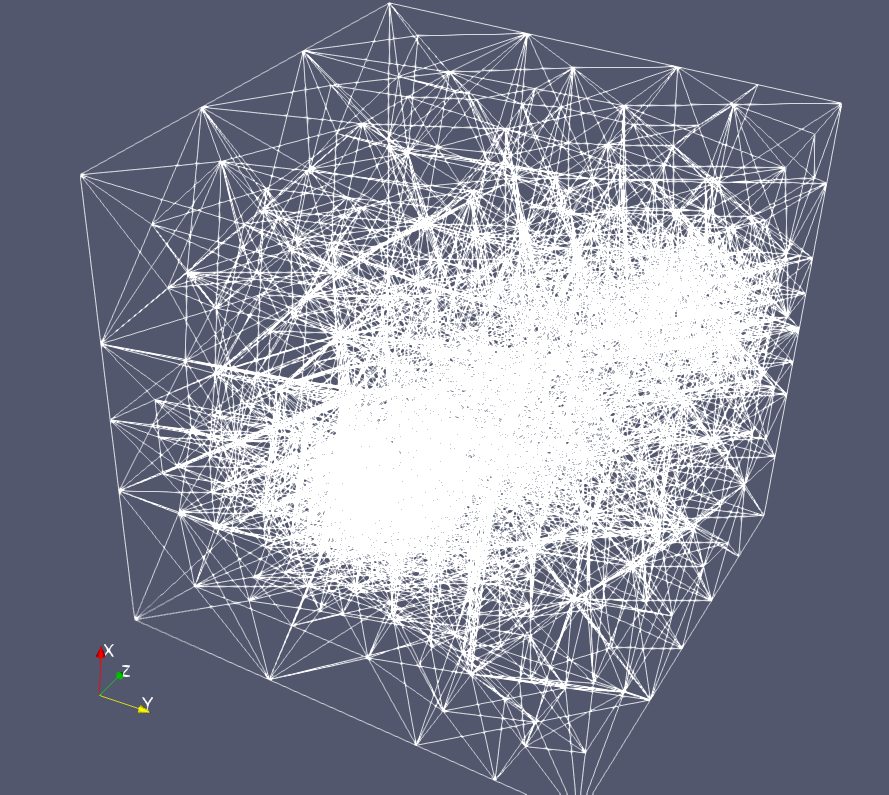} 
  \caption{refinement level 4}
  \end{subfigure}

  \begin{subfigure}[b]{0.46\textwidth}
  \centering
  \includegraphics[width=\textwidth]{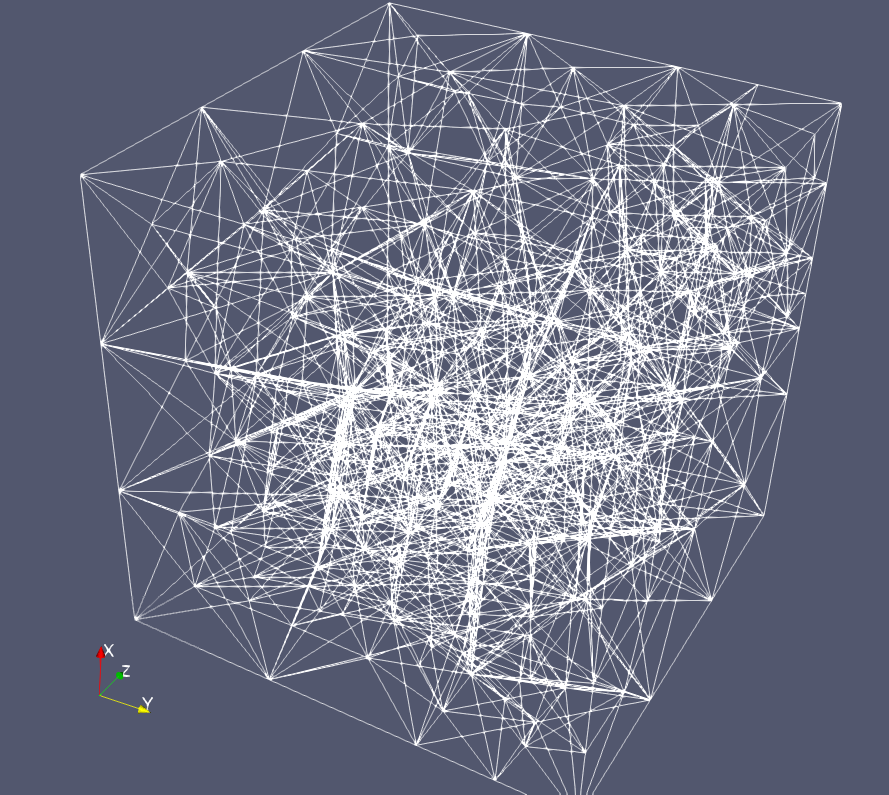} 
  \caption{refinement level 2}
  \end{subfigure}
  ~
  \begin{subfigure}[b]{0.46\textwidth}
  \centering
  \includegraphics[width=\textwidth]{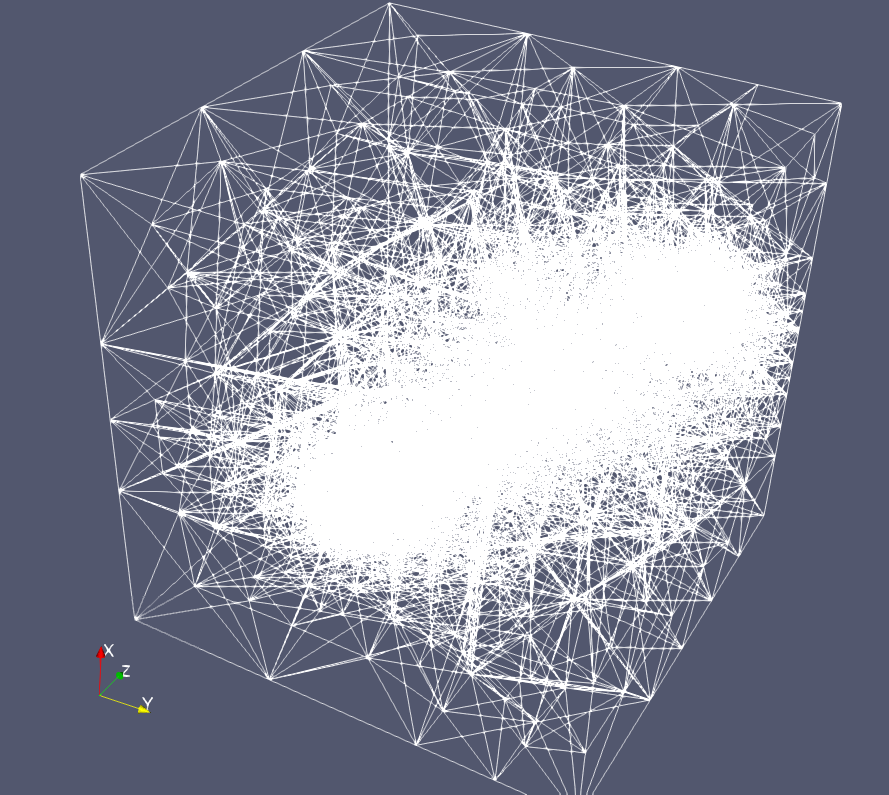} 
  \caption{refinement level 5}
  \end{subfigure}

  \begin{subfigure}[b]{0.46\textwidth}
  \centering
  \includegraphics[width=\textwidth]{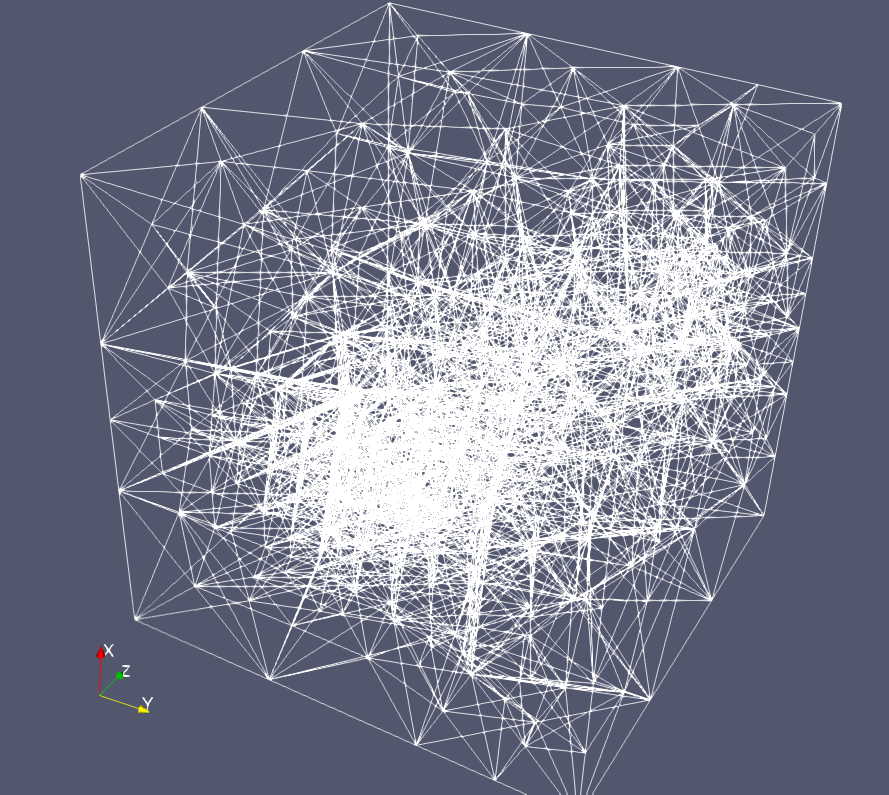} 
  \caption{refinement level 3}
  \end{subfigure}
  ~
  \begin{subfigure}[b]{0.46\textwidth}
  \centering
  \includegraphics[width=\textwidth]{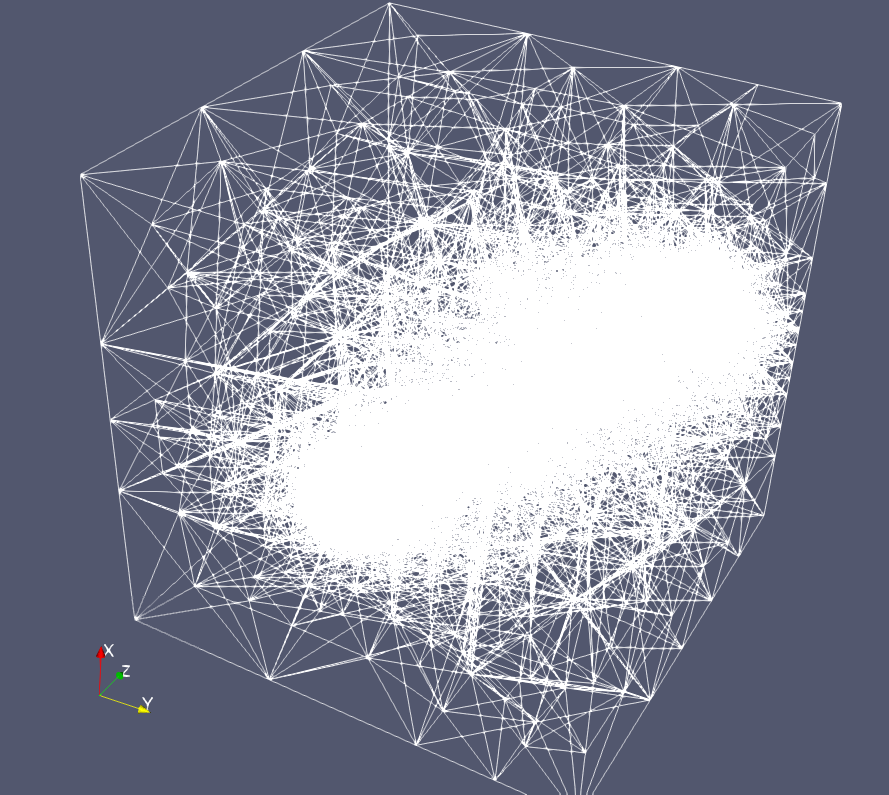} 
  \caption{refinement level 6}
  \end{subfigure}
  
  \caption{Space-time mesh in subsequent refinement steps}
\label{fig:meshes}
\end{figure}

Finally, we include an adaptive finite element method algorithm: Starting with the coarse initial space-time mesh (see Fig.~\ref{fig:init-mesh}),
we iteratively refine the mesh.
As the error indicator, we employ the value of ~$J(\SigU_h, \phi_h)$. As the criterion for selecting elements for refinements, we use the  D\"orfler marking criterion~\cite{st9} with~$\theta = 0.5$.
The resulting convergence of the numerical solution at particular refinement steps is illustrated in Figure \ref{fig:adapt}. 
We can visually estimate the improvement of the quality of the solution.
We also present the convergence of the contour of the solution in Figure \ref{fig:cont}.
We can see an improvement in the quality of the contour as we refine the mesh.
Finally, the sequence of the refined space-time meshes is presented in Figure \ref{fig:meshes}.
The detailed convergence analysis as well as the execution times are listed in Table 2.

\begin{table}[t]
\centering
\begin{tabular}{cccc}
  \hline
  {level} & {DoFs} & {$J(\SigU_h, \phi_h)$} & {solver [{s}]} \\
  \hline
0 &    3,019 & 0.0137 &    0.046  \\
1 &   11,540 & 0.0039 &    0.567  \\
2 &   25,963 & 0.0052 &    4.280  \\
3 &   63,033 & 0.0040 &   14.573  \\
4 &  163,755 & 0.0031 &   31.190  \\
5 &  359,658 & 0.0023 &  101.335 \\
6 &  730,953 & 0.0017 &  279.838  \\
  \hline
\end{tabular}
\caption{Results of the adaptive mesh refinement process}
\end{table}

\section{Conclusions and Future work}\label{Sec:conclusions}

The constraint least-square formulation for the space-time finite element method allows the development of a solver algorithm that requires the inversion of the space-time mass matrix $\color{red}{M}$ and the application of an iterative solver for the Uzawa type of system of equations.
For the space-time advection-dominated diffusion problem, this is possible due to the following algebraic properties: First, $-\color{cadmiumgreen}{L_y^T}\color{purple}{M^{-1}}\color{cadmiumgreen}{L_y}
  -\color{cadmiumgreen}{L_x^T}\color{purple}{M^{-1}}\color{cadmiumgreen}{L_x}+\color{red}{K}=\color{black}{0}$, where $\color{red}{K}$ and $L_x,L_y$ are the space-time matrices defined by (\ref{eq:matrices}).
Second, $(
  \color{blue}{A_x}\color{purple}{M}^{-1}\color{blue}{A_x^T}+\color{blue}{A_y}\color{purple}{M^{-1}}\color{blue}{A^T_y}+
  \color{blue}{A_t}\color{purple}{M}^{-1}\color{blue}{A_t^T})\color{black}{=S}$, where $S$ denotes the space-time stiffness matrix. 
  At some point in the solution process, we deal with $S^{-1}$, which requires that the solution is fixed (e.g. to zero) at some point in the computational domain.
  Finally, we end up with the Uzawa kind of problem 
$ \begin{bmatrix}
  -S & Z \\
  Z^T & 
  0  \\   \end{bmatrix} 
    \begin{bmatrix}
   \lambda \\ \phi 
  \end{bmatrix}
  =
  \begin{bmatrix}
  f \\ 0
  \end{bmatrix},$
where $Z=(\color{cadmiumgreen}{L_y^T} \color{purple}{M^{-1}}\color{blue}{A_y^T}
 + \color{cadmiumgreen}{L_x^T} \color{purple}{M^{-1}}\color{blue}{A_x^T}
  +\color{purple}{M_{\phi\sigma}}\color{purple}{M}^{-1}\color{blue}{A_t^T})$
and
$(\color{blue}{A_t}\color{purple}{M^{-1}M^T_{\phi\sigma}}+\color{blue}{A_x}\color{purple}{M}^{-1}\color{cadmiumgreen}{L_x}
  +\color{blue}{A_y}\color{purple}{M}^{-1}\color{cadmiumgreen}{L_y}
  )S^{-1}(\color{cadmiumgreen}{L_y^T} \color{purple}{M^{-1}}\color{blue}{A_y^T}
 + \color{cadmiumgreen}{L_x^T} \color{purple}{M^{-1}}\color{blue}{A_x^T}
  +\color{purple}{M_{\phi\sigma}}\color{purple}{M}^{-1}\color{blue}{A_t^T})$ is symmetric and positive definite.
This Uzawa kind of system can be solved by an iterative solver, e.g. GMRES algorithm.
We have verified our findings on the IGA discretization and the adaptive finite element method.

\section{Acknowledgement}

Research project supported by the program ``Excellence initiative - research university" for the AGH University of Science and Technology.

\end{document}